\documentclass[twocolumn]{autart}    % Enable this line and disable the 
\usepackage{graphicx}          % Include this line if your 
\usepackage{graphicx}      % include this line if your document contains figures
\usepackage{amsmath,mathtools,amssymb}
\usepackage{amsfonts}
\usepackage{xcolor}
\usepackage{tikz}
\usepackage{pgfplots}
\usetikzlibrary{external}
\usepackage{epstopdf}
\usepackage[font=footnotesize, labelfont={sf,bf}, margin=1cm]{caption}
\usepackage{graphicx, import}
\usepackage{algorithm}
\usepackage{algpseudocode}
\usepackage{booktabs}
\usepackage{caption}
\usepackage{comment}
\usepackage{enumitem}
\usepackage{accents}
\usepackage{hyperref}
\usepackage{scalerel}
\usepackage{amssymb}

\newif\ifcommentandrea
\commentandreafalse
\newcommand{\andrea}[1]{\ifcommentandrea{\textcolor{red}{(andrea: #1)}}\else \fi}

\newif\ifcommentthesis
\commentthesisfalse
\newcommand{\thesis}[1]{\ifcommentthesis{\textcolor{blue}{(fix for thesis: #1)}}\else \fi}

\newcommand{\norm}[1]{\left\lVert#1\right\rVert}

\newcommand{\MR}[1]{\mathrm{#1}}

\makeatletter
\def\BState{\State\hskip-\ALG@thistlm}
\makeatother

\newcommand{\mylabel}[2]{#2\def\@currentlabel{#2}\label{#1}}
\usetikzlibrary{decorations.pathreplacing,calc}

\newcommand{\en}[0]{}
\begin{document}

\begin{frontmatter}

\title{A Lyapunov Function for the Combined System-Optimizer \\Dynamics in Inexact Model Predictive Control} % Title, preferably not more 

\thanks{This research was supported by the German Federal 
Ministry for Economic Affairs and Energy (BMWi) via 
eco4wind (0324125B) and DyConPV (0324166B), by 
DFG via Research Unit FOR 2401, by the EU via ITN-AWESCO
(642 682), and by the Office of Naval Research, ONR grant No. N00014-20-1-2088.}
%\thanks[footnoteinfo]{This paper was not presented at any IFAC 
%	meeting. Corresponding author M.~T.~Cicero. Tel. +XXXIX-VI-mmmxxi. 
%	Fax +XXXIX-VI-mmmxxv.}

\author[ETHZ]{Andrea Zanelli}\ead{zanellia@ethz.ch},    % Add the 
\author[UNC]{Quoc Tran-Dinh}\ead{quoctd@email.unc.edu},    % Add the 
\author[ALUFR]{Moritz Diehl}\ead{moritz.diehl@imtek.uni-freiburg.de}               % e-mail address 
\vspace{0.1cm}
\address[ETHZ]{ETH Zurich, Switzerland}
\address[ALUFR]{University of Freiburg, Germany}
\address[UNC]{The University of North Carolina at Chapel Hill, NC-27599, USA}
                             												% keyword list or with the 
																                        % keyword list or with the 
												                                        % help of the Automatica 
												                                        % keyword wizard
						                                        
\begin{abstract}                         
In this paper, an asymptotic stability proof for a class of methods 
for inexact nonlinear model predictive control is presented. 
General Q-linearly convergent online optimization methods are considered and an asymptotic stability result
is derived for the setting where a limited number of iterations of the optimizer are carried 
out per sampling time. Under the assumption of Lipschitz continuity of the solution, we explicitly construct 
a Lyapunov function for the combined system-optimizer dynamics, 
which shows that asymptotic stability can be obtained if the sampling time is sufficiently short. 
The results constitute an extension to existing attractivity results which hold in the 
simplified setting where inequality constraints are either not present or inactive in the region 
of attraction considered. Moreover, with respect to the established results on robust asymptotic 
stability of suboptimal model predictive control, we develop a framework that takes into account the optimizer's dynamics 
and does not require decrease of the objective function across iterates.
By extending these results, 
the gap between theory and practice of the standard real-time iteration strategy is 
bridged and asymptotic stability for a broader class of methods is guaranteed.
\end{abstract}

\end{frontmatter}

% \textcolor{red}{ 
% \begin{itemize}
%     \item Recursive feasibility: Use soft constraints? Assume (see next point)?
%     \item Solvability of QP? Solution: Theorem~\ref{thm:contraction} + \cite[Lemma 3.5]{TranDinh2012b} 
%     \item Proof of Proposition \ref{prop:bilinear_growth}? 
% \end{itemize}
% }

\section{Introduction}
% \andrea{@Moritz: I am still not sure about the name ``real-time methods''. I think it would be easy for people to make the connection with the real-time iteration, on the other side it might 
% be confusing because of the specific meaning that ``real-time'' has in the context of real-time computing ($\approx$ certifiable) and the generic meaning it has sometimes in the field of 
% model predictive control (real-time optimization $\approx$ fast(-er))}.
\par

Nonlinear model predictive control (NMPC) is an advanced control technique that 
requires the solution of a series of nonlinear programs in order to evaluate an 
implicit control policy. Due to the potentially prohibitive computational burden 
associated with such computations, efficient methods for the solution of the 
underlying optimization problems are of crucial importance. 
For this reason, NMPC was first proposed and applied to systems with slow 
dynamics such as chemical processes in the 70s (see, e.g., \cite{Rawlings2017}).
The interest drawn in both industry and academia has driven in the past decades 
a quick progress in both algorithms and software implementations. At the same time, 
the computational power available on embedded control units has drastically 
increased leading to NMPC gradually 
becoming a viable solution for applications with much shorter sampling times (see, e.g., \cite{Zanelli2019b}). 
\par
In order to mitigate the computational requirements, many applications with high sampling rates rely on approximate 
feedback policies. Among other 
approaches, the real-time iteration (RTI) strategy proposed in 
\cite{Diehl2002} exploits a single iteration of a sequential quadratic 
programming (SQP) algorithm in order to compute an approximate solution of 
the current instance of the nonlinear parametric optimization problem. 
By using this solution to warmstart the SQP algorithm at the next sampling time, 
it is possible to track an optimal solution and eventually converge to it, as the system is 
steered to a steady state \cite{Diehl2005b}. Attractivity proofs for the 
RTI strategy in slightly different settings, and under the assumption that 
inequality constraints are either absent or inactive in a neighborhood of the equilibrium, 
are presented in \cite{Diehl2007b} and \cite{Diehl2005b}. In the same spirit, similar algorithms that rely on 
a limited number of iterations are present in the literature. In \cite{Graichen2010}, a 
general framework that covers methods with linear contraction in the objective 
function values is analyzed and an asymptotic stability proof is provided. 
The recent work in \cite{Liao-McPherson2020} addresses a more general setting 
where an SQP algorithm is used. A proof of local 
input-to-state stability is provided based on the assumption that a sufficiently 
high number of iterations is carried out per sampling time. 
% \andrea{add connection between continuous Lyapunov functions and ISS-Lyapunov functions?}. 
In the convex setting, the works in \cite{Feller2017} and \cite{VanParys2018}
introduce stability results for relaxed barrier \textit{anytime} methods and real-time 
projected gradient methods, respectively. Finally, the works in \cite{Scokaert1999}, 
\cite{Pannocchia2011} and, \cite{Allan2017} analyze conditions under which suboptimal NMPC is 
stabilizing given that a feasible warmstart is available. 
\par
All of the above mentioned methods make use of a common idea. A limited number or, in the limit, a single iteration 
of an optimization algorithm are carried out in order to ``track'' the parametric optimal solution while reducing the computational footprint of the 
method. We will refer to these methods as \textit{real-time methods} in order to make an explicit semantic connection with the 
well known RTI strategy \cite{Diehl2002}.
\par
\begin{figure}
\centering
\includegraphics[scale=0.8]{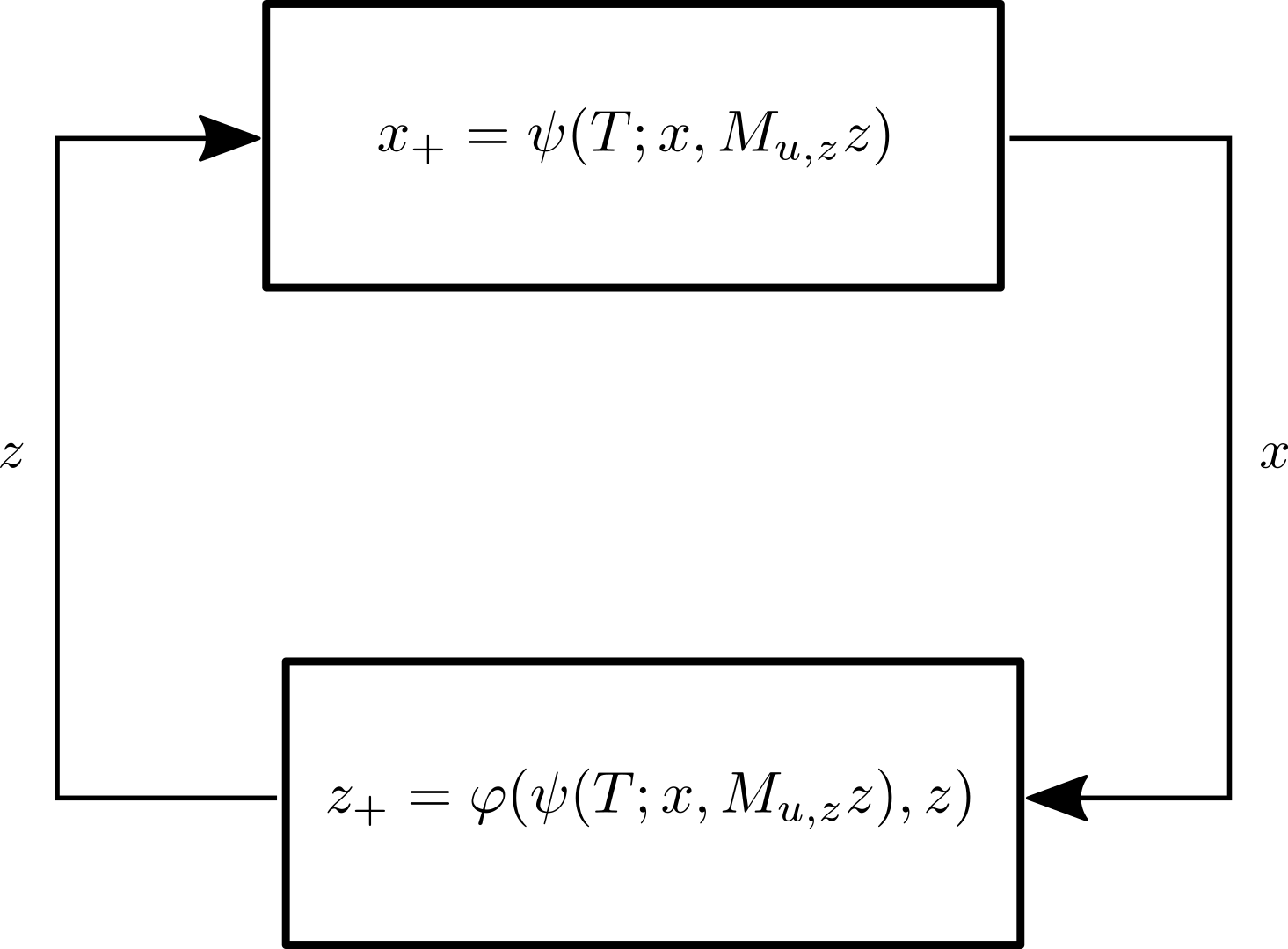}
\caption{Coupled system-optimizer dynamics: when a limited number of iterations 
of the optimization algorithm are carried out in order to obtain an approximate solution that is then 
used to control the system, the system's and the optimizer's dynamics interact with one another.}
\label{fig:system_optimizer}
\end{figure}

Loosely speaking, the main challenge present in real-time approaches lies in the 
fact that the dynamics of the system and the ones of the optimizer interact with each other 
in a non-trivial way, as visualized in Figure \ref{fig:system_optimizer}. 
Although a formal 
definition of the system-optimizer dynamics requires the introduction of several concepts that we 
delay to Section \ref{sec:nominal_so_dyn}, loosely speaking, for a given 
state $x$ and an approximate solution $z$, the system is controlled using the control $u=M_{u,z}z$
as feedback law and, and after the sampling time $T$, it is steered to $x_+ = \psi(T;x, M_{u,z}z)$.
Analogously, the optimizer generates a new approximate solution $z_+ = \varphi(\psi(T;x, M_{u,z}z), z)$. 
We will refer to these coupled dynamics, with state $\xi = (x,z)$, as $\xi_+ = \Phi(T;\xi)$, which will be formally defined in Definition 
\ref{def:sod}.
\color{black}
\subsection{Contribution and Outline}
In this paper, we regard general \textit{real-time methods} and establish asymptotic stability of 
the closed-loop system-optimizer dynamics $\xi_+ = \Phi(T;\xi)$. We assume that the \textit{exact} solution to the underlying optimal control problems yields an 
asymptotically stable closed-loop system and that the iterates of the real-time method contract Q-linearly. Moreover, 
we assume that the primal-dual solution to the underlying parametric nonconvex program is Lipschitz continuous. Under these 
assumptions and the additional assumption that the sampling time $T$ of the closed-loop is sufficiently short, we show asymptotic stability of the system-optimizer dynamics by constructing a Lyapunov function for  $\xi_+ = \Phi(T;\xi)$ and the equilibrium $\xi = 0$ in 
a neighborhood of $\xi = 0$.
\par
In the setting of \cite{Scokaert1999}, 
\cite{Pannocchia2011}, and \cite{Allan2017}, no regularity assumptions 
are required for the optimal solution and optimal value function, which are even allowed to be discontinuous. 
However, a decrease in the objective function is required in order for the optimizer's iterates to be accepted.
This condition is in general difficult to satisfy given that commonly used numerical methods do not generate 
feasible iterates and, for this reason, it is not easy to enforce decrease in the objective function. Although 
robust stability could still be guaranteed by shifting the warmstart (cf. \cite{Allan2017}), the improved iterates might 
be rejected unnecessarily. Moreover, the optimizer's dynamics are completely neglected. With respect to \cite{Scokaert1999}, 
\cite{Pannocchia2011}, and \cite{Allan2017}, we propose in this work an analysis that, although requires stronger assumptions 
on the properties of the optimal solution, incorporates knowledge on the optimizer's dynamics and does not require a decreasing cost. 
\par
Notice that attractivity proofs for the real-time iteration strategy are derived in \cite{Diehl2007b} and \cite{Diehl2005b}
for a simplified setting where inequality constraints are not present or inactive in the entire region of attraction of the closed-loop system.
In this sense, the present paper extends the results in \cite{Diehl2005b} and \cite{Diehl2007b} to a more general setting where active-set changes are allowed within the region of attraction. Moreover, 
asymptotic stability, rather than attractivity, is proved. 
\par
Finally, with respect to \cite{Liao-McPherson2020}, we analyze how the behavior of the system-optimizer dynamics 
is affected by the sampling time, rather than by the number of optimizer's iterations carried out. Additionally, we explicitly construct 
a Lyapunov function for the system-optimizer dynamics. 

\subsection{Notation}
Throughout the paper we denote the Euclidean norm and the $\ell_1$ norm by $\| \cdot \|$ and 
$\| \cdot\|_1$, respectively. We will sometimes write $\| \cdot \|_2$ explicitly, to denote the 
Euclidean norm, when it improves clarity in the derivations. All 
vectors are column vectors and we denote the concatenation of two vectors 
by 
\begin{equation}
    (x,y)\vcentcolon=\begin{bmatrix}x \\ y\end{bmatrix}.
\end{equation}
We denote the derivative (gradient) of any function by 
$\nabla f(x) = \frac{\partial f}{\partial x}(x)^{\top}$ and the Euclidean ball 
of radius $r$ centered at $x$ as $\mathcal{B}(x,r) \vcentcolon=\{ y \,\vcentcolon \,\| x - y \| \leq r\}$.
We use $\mathbb{R}^{m\times n}_{>}$  ($\mathbb{R}^{n}_{>}$) to denote the 
space of strictly positive matrices (vectors) (i.e., the space of matrices (vectors) whose elements are 
real and strictly positive). With a slight abuse of notation we will sometimes write $v > 0$, 
to indicate that all the components of the vector $v$ are strictly positive. Analogously, 
we write $\mathbb{R}^{m\times n}_{\geq}$ ($\mathbb{R}^{n}_{\geq}$) to denote the 
space of matrices (vectors) whose elements are nonnegative and use $v \geq 0$ to indicate that all 
the components of the vector $v$ are nonnegative. We denote 
a vector whose components are all equal to 1 as $\mathbf{1}$ and the identity matrix as $\mathbb{I}$.
Finally, we denote the Minkowski sum $C = \{ c \vcentcolon c = a + b, a \in A, b \in B\}$ of two sets $A$ and $B$ as
$C = A \oplus B$.
\section{System and Optimizer Dynamics}\label{sec:nominal_so_dyn}
In this section, we will define the nominal dynamics that 
the system and optimizer state obey independently from each other. The nominal
dynamics of the closed-loop system are assumed to be such that a Lyapunov function can be constructed if the
exact solution to the underlying discretized optimal control problem is used as feedback law. Similarly, 
we will assume that certain contraction properties hold for the iterates generated by the optimizer if the 
parameter describing the current state of the system is held fixed.

\subsection{System and optimizer dynamics}
In order to study the interaction between the system to be controlled
and the optimizer, we will first formally define their dynamics and describe the 
assumptions required for the stability analysis proposed.
\subsubsection{System dynamics}
\par
The system under control obeys the following sampled-feedback 
closed-loop dynamics:
\par
\begin{defn}[System dynamics]\label{def:ts}
    Let the following differential equation describe the dynamics of the system 
    controlled using a constant input $u_0$: 
    \begin{equation}
        \begin{aligned}
            &\frac{d\psi}{dt}(t;x_0, u_0) = \phi(\psi(t; x_0, u_0), u_0), \\
            &\psi(0; x_0, u_0) = x_0. 
    \end{aligned}
    \end{equation}
    Here, $\psi \vcentcolon \mathbb{R} \times \mathbb{R}^{n_x} \times \mathbb{R}^{n_u} \rightarrow \mathbb{R}^{n_x}$ describes the 
    trajectories of the system, $x_0$ denotes the state of the system at a given sampling instant and 
    $u_0$ the corresponding constant input.
    We will refer to the strictly positive parameter $T > 0$ as the \textit{sampling time} associated with the corresponding discrete-time system 
    \begin{equation}
        x_{\mathrm{next}} = \psi(T; x, u).
    \end{equation}
\end{defn}
We will assume that a slightly tailored type of Lyapunov function
is available for the closed-loop system controlled with a specific policy. 
\begin{assum}\label{ass:lyapunov_stability}
    Let $\bar{u}\vcentcolon\mathbb{R}^{n_x} \rightarrow \mathbb{R}^{n_u}$, and let 
    $V \vcentcolon \mathbb{R}^{n_x} \rightarrow \mathbb{R}$ be a continuous function. 
    Let $\bar{V}$ be a strictly positive constant and define 
    \begin{equation}
        X_{\bar{V}} \vcentcolon = \{ x \vcentcolon V(x) \leq
        \bar{V}\}.
    \end{equation}
    Assume that there exist positive constants $a_1,\, a_2, a_3, T_0$, and
    $q \in [1,\infty)$ such that, for any $x \in X_{\bar{V}}$ and any $T \leq T_0$,
    the following hold:
    \begin{subequations}
        \begin{align}
            a_1 \|x\|^q \leq V(x) &\leq a_2 \|x\|^q, \label{equ:assum_lyap_1} \\ 
            V(\psi(T; x, \bar{u}(x))) - V(x) &\leq -T \cdot a_3 \|x\|^q. \label{equ:assum_lyap_2}
        \end{align}
    \end{subequations}
\end{assum}
Notice that Assumption \ref{ass:lyapunov_stability}, for a fixed $T$ boils down to 
the standard assumption for exponential asymptotic stability (see,
e.g., Theorem 2.21 in \cite{Rawlings2017}). Moreover, 
the dependency on $T$ in \eqref{equ:assum_lyap_2} can be justified, for example, by assuming that a 
continuous-time Lyapunov function $V_{\MR{c}}$ exists such that $\frac{\MR{d}}{\MR{d}t} V_{\MR{c}}(x(t)) \leq -\underline{a} \| x \|^2$, for some 
positive constant $\underline{a}$ and that $V$ is a sufficiently good approximation. 
Moreover, we make the following assumption 
which establishes additional regularity properties of the Lyapunov function.
\begin{assum}\label{ass:qlip}
    Assume that $V^{\frac{1}{q}}$ is Lipschitz continuous over $X_{\bar{V}}$,
    i.e., there exists a constant $\tilde{\mu} > 0$ such that
    \begin{equation}
        |V(x')^{\frac{1}{q}} - V(x)^{\frac{1}{q}}| \leq \tilde{\mu}\| x' - x\|,
    \end{equation}
    $\forall x, x' \in X_{\bar{V}}$.
\end{assum}

\begin{rem}
    Notice that a sufficient condition for Assumption \ref{ass:qlip} to hold is that 
    $V$ is Lipschitz continuous over $X_{\bar{V}}$ and that $V^{\frac{1}{q}}$ is 
    Lipschitz continuous at $x=0$. These conditions are satisfied, for example, by Lyapunov functions 
    which are twice continuously differentiable at the origin if $q=2$ or 
    simply Lipschitz continuous at the origin if $q=1$.
\end{rem}
The following proposition provides asymptotic stability of the
closed-loop system obtained using the feedback policy
$\bar{u}$.
\begin{prop}[Lyapunov stability]\label{pro:conv_stability}
    Let Assumption \ref{ass:lyapunov_stability} hold. Then, for any $T \leq
    T_0$, the origin is an exponentially asymptotically stable equilibrium 
    with region of attraction $X_{\bar{V}}$ for the closed-loop
    system $x_{\MR{next}} = \psi(T;x, \bar{u}(x))$.
    \begin{pf}
        Due to Assumption \ref{ass:lyapunov_stability}, the function $V$ is a valid Lyapunov function for the
        closed-loop dynamics for any $T \leq T_0$. \hfill $\qed$
    \end{pf}
\end{prop}
The Lyapunov function defined in Assumption
\ref{ass:lyapunov_stability} guarantees that, if the ideal policy
$\bar{u}$ is employed, the resulting closed-loop system is
(exponentially) asymptotically stable. In the following, we define the
dynamics of the optimizer used to
numerically compute approximations of $\bar{u}(x)$ for a given state $x$.
\subsubsection{Optimizer dynamics} 
We will assume that we dispose of a numerical method that defines what
we will call the \textit{optimizer} (or more generally a
\textit{solver}) that, for a given $x$, can compute a vector
$\bar{z}(x)$ from which we can compute $\bar{u}(x)$ through a linear
map.
\begin{assum}\label{ass:policy_solver}
    Assume that there exists a function $\bar{z}
    \vcentcolon X_{\bar{V}} \rightarrow \mathbb{R}^{n_z}$ and a matrix
    $M_{u,z}$ such that, for any $x \in X_{\bar{V}}$, the following
    holds:
    \begin{equation}
        \bar{u}(x) = M_{u,z}\bar{z}(x).
    \end{equation}
    For simplicity of notation, we will assume further that $\|
    M_{u,z} \| = 1$.
\end{assum}
\andrea{add remark on localization?}
\begin{defn}[Optimizer dynamics]\label{def:od}
    Let the following discrete-time system describe the dynamics of the optimizer
    \begin{equation}\label{equ:od}
        z_+ = \varphi(\psi(T; x, M_{u,z}z), z), \\
    \end{equation}
    where $\varphi \vcentcolon \mathbb{R}^{n_x} \times \mathbb{R}^{n_z} \rightarrow \mathbb{R}^{n_x}$ 
    and $z$ represents the state of the optimizer.
\end{defn}

\begin{rem}
    Notice that the optimizer dynamics \eqref{equ:od} make use of the current approximate solution $z$
    and a forward-simulated state $x_+ = \psi(T; x, M_{u,z}z)$. This setting corresponds, for example,
    to the case where a real-time iteration is carried out by solving a QP where the parameter $x_+ = \psi(T; x, M_{u,z}z)$
    is used as current state of the system. This amounts to assuming that either a perfect prediction $x_+$ of the system's state is 
    available ahead of time, or that instantaneous feedback can be delivered to the system. In both cases,
    small perturbations introduced by either model mismatch or feedback delay could be introduced explicitly. 
    This goes however beyond the scope of the present work.
\end{rem}
In order to leverage a certain type of contraction estimates, we will make the two following assumptions.
\par
\andrea{need to treat boundaries here! done!}
\thesis{need to treat boundaries here!}
\begin{assum}[Lipschitz continuity]\label{ass:Lip_cont_z} Assume that
    there exist positive constants
    $\hat{r}_x$ and $\sigma$ such that, for any $x \in X_{\bar{V}}$ and any $x' \in \mathcal{B}_{\en}(x, \hat{r}_x)$, the following
    holds
    \begin{equation}
        \| \bar{z}(x') - \bar{z}(x) \| \leq \sigma \| x' - x\|. 
    \end{equation}
    Moreover, we assume that $\bar{z}(0) = 0$.
\end{assum}

\begin{assum}[Contraction]\label{ass:contraction}
    There exist positive constants $\hat{r}_z > 0$ and $\hat{\kappa} < 1$ such that, for any $x\in X_{\bar{V}}$ 
    and any $z \in \mathcal{B}(\bar{z}(x), \hat{r}_{z})$, the
    optimizer dynamics produce $z_+$ such that
    \begin{equation}\label{equ:contraction}
        \norm{z_+ - \bar{z}(x)} \leq \hat{\kappa} \norm{z - \bar{z}(x)}.
    \end{equation}
\end{assum}
The following lemma provides a way of quantifying the perturbation to the nominal
contraction \eqref{equ:contraction} due to changes in the value of
$x$ across iterations of the optimizer.
\begin{lem}\label{lem:contraction}
    Let Assumptions 
    \ref{ass:Lip_cont_z} and \ref{ass:contraction} hold. Then there exist strictly positive 
    constants $r_z \leq \hat{r}_z$ and $r_x \leq \hat{r}_x$ such that, 
    for any $x$ in $X_{\bar{V}}$, any $z$ in $\mathcal{B}(\bar{z}, r_z)$, and 
    any $x'$ in $\mathcal{B}_{\en}(x, r_x)$, we have that
    \begin{align}\label{equ:track_contraction}
     \begin{split}
         \norm{z_+ - \bar{z}(x')} & \leq  \hat{\kappa} \norm{z -
         \bar{z}(x)} +  \sigma \hat{\kappa} \norm{x' - x}_{\en}.
        \end{split} 
    \end{align}
    \begin{pf}
        By choosing $r_x$ and $r_z$ such that they satisfy $0 < r_x \leq \hat{r}_x$ and 
        $0 < r_z \leq \hat{r}_z - \sigma r_x$, any $x$ in $X_{\bar{V}}$,  any $z \in \mathcal{B}(\bar{z}, r_z)$, and any 
        $x' \in \mathcal{B}(x,r_x)$, we have that $\|z - \bar{z}(x')\| \leq \hat{r}_z$. Hence, we can apply
        the contraction from Assumption \ref{ass:contraction} together with Lipschitz continuity of $\bar{z}$ from Assumption \ref{ass:Lip_cont_z} 
        in order to write
        \begin{equation}
            \begin{aligned}
                &\| z_+ - \bar{z}(x')\| &&\leq \hat{\kappa} \| z - \bar{z}(x')\| \\
                &                       &&= \hat{\kappa}\| z - \bar{z}(x) + \bar{z}(x) - \bar{z}(x')\| \\
                &                       &&\leq \hat{\kappa}\| z - \bar{z}(x) \| + \hat{\kappa} \| \bar{z}(x) - \bar{z}(x')\| \\
                &                       &&\leq \hat{\kappa}\| z - \bar{z}(x) \| + \sigma \hat{\kappa} \| x - x'\|, \\
            \end{aligned}
        \end{equation}
        which concludes the proof. \hfill$\qed$
    \end{pf}
\end{lem}
\subsection{Discussion of fundamental assumptions}
Notice that the setting formalized by Assumptions \ref{ass:policy_solver}, \ref{ass:Lip_cont_z}, and
\ref{ass:contraction} does not require that $V$ is
the optimal value function of a discretized optimal control problem nor of an
optimization problem in general. We require instead that it is a
Lyapunov function with some additional properties according to 
Assumptions \ref{ass:lyapunov_stability} and \ref{ass:qlip}. 
Similarly, $\bar{z}$ needs not be the primal-dual solution to an
optimization problem. We require instead that it is associated
with the policy $\bar{u}$, i.e., for any $x \in X_{\bar{V}}$, 
$\bar{u}(x) = M_{u,z}\bar{z}(x)$, 
that it is Lipschitz continuous and that the optimizer (or
``solver'' in general) can generate Q-linearly contracting
iterates that converge to $\bar{z}(x)$.
\par
% \begin{rem}
    However, in order to make a more concrete connection with a classical setting, in NMPC
    we can often assume that $V$ is the optimal value
    function of a discretized version of a continuous-time optimal control problem. 
    In this case, we can refer to the resulting
    parametric nonlinear program of the following form
\begin{equation}\label{eq:nmpc}
{\!\!\!\!\!}{\!\!}\begin{aligned}
&\underset{\begin{subarray}{c}
s_0, \dots, s_N \\
u_0, \dots, u_{N-1}
\end{subarray}}{\min}	    &&\sum_{i=0}^{N-1} l(s_i, u_i) + m(s_N)\\ 
&\,\,\,\quad \text{s.t.}    &&s_0 - x = 0, \\
& 						    &&\psi_{\text{d}}(s_i,u_i)  - s_{i+1} = 0, \,\,\, i = 0,\dots, N-1,\\
& 						    &&\pi(s_i, u_i) \leq 0, \quad \quad \quad \,\,\,\,\,\, i = 0,\dots, N-1,\\
& 						    &&\pi_N(s_N) \leq 0, 
\end{aligned}{\!\!\!}
\end{equation}
where $s_i \in \mathbb{R}^{n_s}$, for $i = 0,\dots,N$, 
and $u_i \in \mathbb{R}^{n_u}$, for $i = 0,\dots,N-1$, 
describe the predicted states and inputs of
the system to be controlled, respectively. The functions 
$l \vcentcolon \mathbb{R}^{n_s} \times \mathbb{R}^{n_u} \rightarrow \mathbb{R}$ and  
$m \vcentcolon \mathbb{R}^{n_s} \rightarrow \mathbb{R}$ describe the stage and 
terminal cost, respectively, and $\psi_{\text{d}} \vcentcolon \mathbb{R}^{n_s} \times \mathbb{R}^{n_u} \rightarrow \mathbb{R}^{n_s}$, 
$\pi \vcentcolon \mathbb{R}^{n_s} \times \mathbb{R}^{n_u} \rightarrow \mathbb{R}^{n_{\pi}}$, and
$\pi_N \vcentcolon \mathbb{R}^{n_s} \rightarrow \mathbb{R}^{n_{\pi_{\scaleto{N}{2pt}}}}$ describe the dynamics of the system, 
the stage, and terminal constraints, respectively.

The first-order optimality condition associated with \eqref{eq:nmpc} can be represented 
by a generalized equation (see, e.g., \cite{Robinson1980}), i.e.,
\begin{equation}
    0 \in F(z,x) + \mathcal{N}_K(z)
\end{equation}
where $z \in \mathbb{R}^{n_z}$ denotes the primal-dual solution, 
$F \vcentcolon \mathbb{R}^{n_z}\rightarrow \mathbb{R}^{n_z}$ is a vector-valued map 
and $\mathcal{N}_K \vcentcolon \mathbb{R}^{n_z} \rightarrow \mathbb{R}^{n_z}$ is a set-valued map denoting the normal cone to a convex set $K \subseteq \mathbb{R}^{n_z}$.
    Under proper
    regularity assumptions, we can then obtain that a single-valued
    localization $\bar{z}$ of the solution map of the generalized equation must
    exist. Under the same assumptions, we can usually prove Lipschitz
    continuity of $\bar{z}$ and Q-linear contraction of, for example,
    Newton-type iterations (cf. \cite{Josephy1979}),
    hence satisfying Assumptions \ref{ass:Lip_cont_z} and
    \ref{ass:contraction}. Notice that existence of a Lipschitz continuous single-valued localization of the solution 
    map, does not require the solution map itself to be single-valued or Lipschitz continuous.
    \par
    Similarly, using the standard argumentation for NMPC (see, e.g., \cite{Rawlings2017}), 
    we obtain that the feedback policy associated
    with the global solution to the underlying nonlinear programs is
    stabilizing and that the associated optimal value function is a Lyapunov
    function. In this way, Assumptions \ref{ass:lyapunov_stability},
    \ref{ass:qlip}, and \ref{ass:policy_solver} are satisfied.
    Note that, in this context, an aspect that remains somewhat unresolved is the fact that the
    standard argumentation for the stability analysis of NMPC requires
    that the global optimal solution is found by the optimizer. Hence, in
    order to be able to identify $V$ with the optimal
    value function it is necessary to assume that the single-valued localization $\bar{z}$ attains the
    global minimum for all $x \in X_{\bar{V}}$.
% \end{rem}
\subsection{Combined system-optimizer dynamics}\label{sec:gen_rti}
Proposition \ref{pro:conv_stability} and Lemma~\ref{lem:contraction} provide 
key properties of the system and optimizer dynamics, respectively. In this section, 
we analyze the interaction between these two dynamical systems and how
these properties are affected by such an interplay. 
To this end, let us define the following coupled system-optimizer dynamics.
\par
\begin{defn}[System-optimizer dynamics]\label{def:sod}
    Let the following discrete-time system describe the coupled system-optimizer dynamics:
    \begin{equation}\label{equ:so_dyn}
        \begin{aligned}
            &&x_+ &= \psi(T; x, M_{u,z}z), \\
            &&z_+ &= \varphi(\psi(T; x, M_{u,z}z), z) \\
        \end{aligned}
    \end{equation}
    or, in compact form
    \begin{equation}\label{equ:so_dyn_compact}
        \xi_+ = \Phi(T;\xi),
    \end{equation}
    where $\xi \vcentcolon = (x,z)$ and $\Phi \vcentcolon \mathbb{R} \times \mathbb{R}^{n_x+n_z} \rightarrow \mathbb{R}^{n_x + n_z}$.
\end{defn}
Our ultimate goal is to prove that, for a sufficiently short sampling
time $T$, the origin is a locally asymptotically stable equilibrium for
\eqref{equ:so_dyn_compact}. In order to describe the interaction between the 
two underlying subsystems, we analyze how the Lyapunov decrease \eqref{equ:assum_lyap_2} and the error contraction \eqref{equ:track_contraction}
 are affected.
\subsubsection{Error contraction perturbation}

\par
In the following, we will specialize the result from Lemma \ref{lem:contraction} to
the context where $\| x'-x\|$ is determined by the evolution of the
system to be controlled under the effect of the approximate policy. To this end, we make 
a general assumption on the behavior of the closed-loop system for a bounded value of the numerical error.

\begin{assum}[Lipschitz system dynamics]\label{ass:state_evolution_bound}
    Assume that $\phi(0,0) = 0$ and that positive finite constants $L_{\phi,x}$, $L_{\phi,u}$, and $\rho$ exist 
such that, for all $x',x \in X_{\bar{V}} \oplus \mathcal{B}(0,\rho)$, 
    all $u'=M_{u,z}z', \,u=M_{u,z}z$, with 
    $z', z\in \mathcal{B}(\bar{z}(x), r_z)$, the following holds:
    \begin{equation}
        \| \phi(x', u') - \phi(x,u) \| \leq L_{\phi,x}\| x' - x \| + L_{\phi,u} \| u'-u \|.   
    \end{equation}
\end{assum}
The following propositions establish bounds on the rate at which the
state can change for given $x$ and $z$.
\begin{prop}\label{pro:g_lemma}
    Let Assumption \ref{ass:state_evolution_bound} hold. Then, there exist 
    a positive finite constant $T_1 > 0$ such that
    \andrea{unify norm and ball notation!}
    for all $x \in X_{\bar{V}}$, all $z \in \mathcal{B}(\bar{z}(x), r_z)$, and 
    any $T \leq T_1$, the following holds:
    \begin{equation}
        \begin{aligned}
            &&\| \psi(T; x, M_{u,z}&z) - x\| \leq \\
            && & T \cdot (L_{\psi,x}\| x \| + L_{\psi,u} \| M_{u,z}z \|),   
    \end{aligned}
    \end{equation}
    where $L_{\psi,x}\vcentcolon=e^{L_{\phi,x}T_1}L_{\phi,x}$ and 
    $L_{\psi_,u}\vcentcolon=e^{L_{\phi,x}T_1}L_{\phi,u}$.
    Moreover, for all $x \in X_{\bar{V}}$, all $u'=M_{u,z}z', \,u=M_{u,z}z$, such that 
    $z', z\in \mathcal{B}(\bar{z}(x), r_z)$, and any $T \leq T_1$, the following 
    holds:
    \begin{equation}
        \| \psi(T; x, u') - \psi(T; x, u)\| \leq T \cdot L_{\psi,u} \| u' - u \|.   
    \end{equation}
    \begin{pf}
        See Appendix \ref{app:g_lemma}. \hfill$\qed$
    \end{pf}
\end{prop}
\color{black}
\begin{prop}\label{pro:state_evolution_bound}
    Let Assumptions \ref{ass:Lip_cont_z} and \ref{ass:state_evolution_bound} hold.
    Then, there exist positive finite constants $\eta, \theta$, and $T_2 > 0$, such that
    for any $x \in X_{\bar{V}}$, any $z \in \mathcal{B}(\bar{z}(x), r_z)$, 
    and any $T\leq T_2$, the following holds:
	\begin{equation}\label{equ:state_evolution_bound}
        \|\psi(T; x, u) - x \| \leq T \cdot (\eta\| x \| + \theta \| z - \bar{z}(x)\|).
	\end{equation}
    \begin{pf}
    Define $\eta \vcentcolon = L_{\psi,x} + L_{\psi,u}\sigma$ and $\theta \vcentcolon = L_{\psi,u}$.
    Due to Proposition~\ref{pro:g_lemma} we have that 
    for any $x \in X_{\bar{V}}$, any $z \in \mathcal{B}(\bar{z}(x), r_z)$,
    and any $T\leq T_1$, the following holds:
    \begin{equation*}
        \| \psi(T; x, M_{u,z}z) - x \| \leq T \cdot \left(L_{\psi,x}\| x \| + L_{\psi,u} \| M_{u,z}z \|\right).   
    \end{equation*}
    Hence, due to Assumption \ref{ass:Lip_cont_z}, we can write 
    \begin{equation*}
        \norm{M_{u,z}z} \leq \norm{\bar{z}(x)} + \norm{z - \bar{z}(x)} \leq \sigma \norm{x} + \norm{z - \bar{z}(x)}
    \end{equation*}
    such that the following holds:
    \begin{equation}
        \begin{aligned}
            &&\norm{\psi(T; x,u) - x} \leq &\,\,T(L_{\psi,x} + L_{\psi,u}\sigma) \norm{x} \\
            && & + TL_{\psi,u}\norm{z - \bar{z}(x)}. 
        \end{aligned}
    \end{equation}
    Finally, defining $T_2 \vcentcolon= \min \{ T_0, T_1\}$ completes the proof. \hfill$\qed$
\end{pf}
\end{prop}

Using the bound from Proposition \ref{pro:state_evolution_bound}
together with the contraction from Lemma \ref{lem:contraction} we
obtain the following perturbed contraction.
\begin{prop}\label{prop:contraction_x}
    Let Assumptions \ref{ass:lyapunov_stability}, \ref{ass:Lip_cont_z}, \ref{ass:contraction}, and \ref{ass:state_evolution_bound} 
hold. Moreover, define 
\begin{equation}\label{equ:Ts_contraction}
T_3' \vcentcolon = \min \Bigg\{\frac{r_x}{\eta r_{\bar{V}} + \theta r_z}, \frac{r_z(1 - \hat{\kappa})}{\sigma \hat{\kappa}(\theta r_z + \eta r_{\bar{V}})} \Bigg\},
\end{equation}
where $r_{\bar{V}} \vcentcolon = \left(\frac{\bar{V}}{a_1}\right)^{\frac{1}{q}}$.

Then, for any $x,z$ such that $x \in X_{\bar{V}}$ and $\| z - \bar{z}(x)\| \leq r_z$ and any $T \leq T_3 \vcentcolon= \min \{ T_3', T_2 \}$, the following holds:
\begin{equation}\label{equ:contr_z}
\| z_+ - \bar{z}(x_+)\|  \leq \, \kappa\| z - \bar{z}(x)\| + T \gamma \| x \|, 
\end{equation}
where
\begin{equation}
    \kappa \vcentcolon = \hat{\kappa}(1+T_3 \sigma\theta) < 1 \quad \text{and} \quad \gamma \vcentcolon = \sigma\hat{\kappa}\eta.
\end{equation} 
Moreover, we also have $\| z_+ - \bar{z}(x_+) \| \leq r_z$.
\begin{pf}
    Due to Assumptions \ref{ass:lyapunov_stability} and \ref{ass:state_evolution_bound}, given that $\| z - \bar{z}(x) \| \leq r_z$ and $T \leq T_3 \leq T_2$, we have $\| x_+ - x \| \leq r_x$ for all $x \in X_{\bar{V}}$ 
    if we additionally require that 
    \begin{equation}
        T \leq \frac{r_x}{\eta r_{\bar{V}} + \theta r_z}, 
    \end{equation}
    which provides the first term in the definition of $T_3'$.
    Hence, for all $T \leq T_3$, we can apply the contraction from 
    Lemma \ref{lem:contraction}: 
    \begin{align}
     {\!\!\!\!}\begin{split}
         \norm{z_+ - \bar{z}(x_+)}_{\en} & \leq  \hat{\kappa} \norm{z - \bar{z}(x)}_{\en} +  \sigma \hat{\kappa} \norm{x_+ - x}_{\en}
        \end{split} {\!\!\!\!}
    \end{align}
    and applying the inequality from Proposition \ref{pro:state_evolution_bound}, we obtain
    \begin{equation}
    \| z_+ - \bar{z}(x_+) \|  \leq \, \kappa\| z - \bar{z}(x)\| + T \gamma \| x \|, 
    \end{equation}
    where
    \begin{equation}
        \kappa \vcentcolon = \hat{\kappa}(1+T_3 \sigma\theta) \quad \text{and} \quad \gamma \vcentcolon = \sigma\hat{\kappa}\eta.
    \end{equation} 
    \par
    From this last inequality we see that in order to guarantee that $\| z_+ - \bar{z}(x_+) \| \leq r_z$, we
    must impose that
    \begin{equation}
        T \leq \frac{r_z(1 - \hat{\kappa})}{\sigma \hat{\kappa}(\theta r_z + \eta r_{\bar{V}})},
    \end{equation}
    which provides the second term in the definition of $T_3'$.

    Finally, since, for any $r_{\bar{V}} > 0$, the following holds
    \begin{equation}
        \frac{1-\hat{\kappa}}{\sigma\hat{\kappa}\theta} > \frac{r_z(1 - \hat{\kappa})}{\sigma \hat{\kappa}(\theta r_z + \eta r_{\bar{V}})},    
    \end{equation}
    we have that $\kappa < 1$ for any $T \leq T_3'$. \hfill$\qed$
\end{pf}
\end{prop}
Proposition \ref{prop:contraction_x} shows that under suitable assumptions and, in particular, for any $T \leq T_3$, we can guarantee that the numerical error does not increase after one iteration of the optimizer.
\par
In the following, we will make similar considerations for the behavior of $V(x)$ across iterations.
\subsubsection{Lyapunov decrease perturbation}
Let us analyze the impact of the approximate 
feedback policy $M_{u,z}z$ on the nominal Lyapunov contraction.
Throughout the rest of the section, for the sake of notational simplicity, we will make use of the following shorthand:
\begin{equation}\label{equ:short2}
    V_+(T;x,z) \vcentcolon= V(\psi(T; x, M_{u,z}z))
\end{equation}
to denote the value taken by the optimal cost at the state reached applying the suboptimal control action $M_{u,z}z$ starting 
from $x$. Similarly, we introduce 
\begin{equation}\label{equ:short3}
    E(x,z)\vcentcolon = \|z - \bar{z}(x)\|
\end{equation}
and 
\begin{equation}\label{equ:short4}
    \begin{aligned}
        &&E_+&(T;x,z) \vcentcolon= \\
        && &\| \varphi(\psi(T; x, M_{u,z}z), z) - \bar{z}(\psi(T; x, M_{u,z}z))\|
    \end{aligned}
\end{equation}
to denote the numerical error attained at the ``current'' and at the next iteration of 
the optimizer, respectively, where the error is computed with respect to the exact solution associated 
with the ``current'' and next state of the system. 

% \begin{rem}
%     Notice that $\mu$ is proportional to the sampling time $T$.
% \end{rem}
% \begin{assum}\label{ass:invariance_gen}
%     Assume that the initial iterate $z^0$ is initialized such that 
%     \begin{equation}
%         e^0 = \| z^0 - \bar{z}^0 \| \leq \tilde{r}_z, 
%     \end{equation}
%     where 
%     \begin{equation}
%         \tilde{r}_z \vcentcolon = \min \bigg\{ r_z, \frac{\bar{a}\bar{V}}{\mu}\bigg\},
%     \end{equation}
%     that $x^0$ is in $X_{\bar{V}}$ and that the sampling time $T$ satisfies the following bound: 
%     \begin{equation}\label{equ:ts_tilde_r_z}
%         T \leq \frac{(1-\kappa)\tilde{r}_z \sqrt{a_1}}{\sqrt{\bar{V}}\gamma}.
%     \end{equation}
%     % where
%     % \begin{equation}
%     %     \bar{\gamma} \vcentcolon = (\sigma \hat{\kappa} L_{\psi,u} + \sigma^2 \hat{\kappa} L_{\psi,x}).
%     % \end{equation}
% \end{assum}
% \begin{rem}
%     As noticed before for Assumption \ref{ass:Ts_{\MR{c}}ontraction}, in this case too 
%     the bounds on $T$ can be seen as bounds on the contraction rate $\hat{\kappa}$ for a given $T$ instead.
% \end{rem}
\begin{prop}\label{pro:inexact_lyapunov_contraction}
    Let Assumptions \ref{ass:lyapunov_stability}, \ref{ass:qlip},
    \ref{ass:policy_solver}, \ref{ass:Lip_cont_z}, \ref{ass:contraction}, and \ref{ass:state_evolution_bound} hold. 
    Then, there exists a finite positive constant $\mu$ such that, for any $z \in \mathcal{B}(\bar{z}(x), r_z)$,
    any $x \in X_{\bar{V}}$, and any $T \leq T_1$, the following holds: 
    \begin{equation}
        V_+(T;x,z)  \leq (1 - T \bar{a}) V(x) + T\mu E(x,z), 
    \end{equation}
    where $\bar{a} \vcentcolon=\frac{a_3}{a_2}$.
    \begin{pf}
    Assumption \ref{ass:lyapunov_stability} implies that, for any $x \in X_{\bar{V}}$
    and any $T \leq T_0$ the following holds:
    \begin{equation}
        \begin{aligned}
            &&V(\psi(T; x, M_{u,z}\bar{z}(x))) &\leq V(x) - T a_3 \| x \|^{q} \\
            &&           &\leq V(x) - T \frac{a_3}{a_2} V(x) \\ 
            &&           &=(1 - T \bar{a}) V(x). 
        \end{aligned}
    \end{equation}
    Due to Assumption \ref{ass:qlip}, for any $x', x \in
    X_{\bar{V}}$, we can write 
    \begin{equation*}
        \begin{aligned}
            &&\!\!\!\!\!\!|V(x') - V(x)| &\leq | (V(x')^{\frac{1}{q}})^q - (V(x)^{\frac{1}{q}})^q| \\
            && & = | (V(x')^{\frac{1}{q}} + V(x)^{\frac{1}{q}})(V(x')^{\frac{1}{q}} - V(x)^{\frac{1}{q}})| \\
            && & \leq 2 \bar{V}^{\frac{1}{q}} \tilde{\mu} \| x' - x \|
    \end{aligned}
    \end{equation*}
    and defining $L_V \vcentcolon = 2 \bar{V}^{\frac{1}{q}} \tilde{\mu}$ we obtain 
    \begin{equation}
    |V(x') - V(x)| \leq L_V \| x' - x\|.
    \end{equation}
    Together with Proposition \ref{pro:g_lemma}, this implies that, for any $z \in \mathcal{B}(\bar{z}(x), r_z)$, any $x \in X_{\bar{V}}$,
    and any $T \leq T_1$, we can write:
    \begin{equation*}
        \begin{aligned}
            &&V(\psi(T&;x, M_{u,z}z)) - V(\psi(T;x, M_{u,z} \bar{z}(x))) \\
            && &\leq|V(\psi(T; x, M_{u,z}z)) - V(\psi(T; x, M_{u,z}\bar{z}(x)))| \\
            && &\leq L_V\|\psi(T; x, M_{u,z}z) - \psi(T; x, M_{u,z}\bar{z}(x))\| \\
            && &\leq T L_{\psi,u}L_V\|z - \bar{z}(x)\|, \\
        \end{aligned}
    \end{equation*}
    which implies
    \begin{equation}
        \begin{aligned}
            && V_+(T;x,z) & \leq (1 - T \bar{a}) V(x) + T \mu E(x,z), \\
        \end{aligned}
    \end{equation}
    where $\mu \vcentcolon = L_V L_{\psi,u}$. \hfill$\qed$
    \end{pf}
\end{prop}
% \begin{rem}
%     Notice that $\mu$ is proportional to the sampling time $T$.
% \end{rem}
Using Proposition \ref{pro:inexact_lyapunov_contraction} we can formulate Lemma \ref{lem:invariance} below which 
establishes positive invariance of the following set for the system-optimizer dynamics \eqref{equ:so_dyn}.
\begin{defn}[Invariant set]
    Define the following set:
    \begin{equation*}
        \Sigma \vcentcolon = \{(x, z) \in \mathbb{R}^{n_x + n_z} \vcentcolon V(x) \leq \bar{V}, \| z - \bar{z}(x)\| \leq \tilde{r}_z\},
    \end{equation*}
    where
    \begin{equation}\label{equ:tilde_r_z_def}
        \tilde{r}_z \vcentcolon = \min \bigg\{ r_z, \frac{\bar{a}\bar{V}}{\mu}\bigg\}.
    \end{equation}
\end{defn}
\begin{lem}[Invariance of $\Sigma$]\label{lem:invariance}
    Let Assumptions \ref{ass:lyapunov_stability}, \ref{ass:qlip},
    \ref{ass:policy_solver}, \ref{ass:Lip_cont_z}, \ref{ass:contraction}, 
    and \ref{ass:state_evolution_bound} hold. Define  
    \begin{equation}\label{equ:ts_tilde_r_z}
        T_4' \vcentcolon = \frac{(1-\kappa)\tilde{r}_z a_1^{\frac{1}{q}}}{\bar{V}^{\frac{1}{q}}\gamma}.
    \end{equation}
    Then, for any $\xi \in \Sigma$ and any $T \leq T_4 \vcentcolon= \min \{ T_4', T_3\}$, 
    it holds that $\xi_+ \in \Sigma$. Moreover, the following coupled inequalities hold:
    \begin{equation}\label{equ:coupled_contraction}
        \begin{aligned}
            && V_+(T;x,z) & \leq (1 - T \bar{a}) V(x) + T \mu E(x,z), \\
            && E_+(T;x,z) & \leq  T \hat{\gamma} V(x)^{\frac{1}{q}}  + \kappa  E(x,z), \\
        \end{aligned}
    \end{equation}
    where $\hat{\gamma} \vcentcolon=\gamma / a_1^{\frac{1}{q}}$.
    \begin{pf}
        Given that $E(x,z) \leq \tilde{r}_z \leq r_z$ and $x \in X_{\bar{V}}$, 
        we can apply the contraction 
        from Proposition \ref{pro:inexact_lyapunov_contraction}, such that 
    \begin{equation}
        V_+(T;x,z)  \leq (1- T \bar{a}) V(x) + T \mu E(x,z), 
    \end{equation}
    holds. Moreover, due to the definition of $\tilde{r}_z$ in \eqref{equ:tilde_r_z_def}, we have that
    $V_+(T;x,z) \leq \bar{V}$ since
    \begin{equation}
        \begin{aligned}
            &&V_+(T;x,z)  &\leq (1- T \bar{a}) V(x) + T \mu E(x,z)\\
            &&            &\leq (1- T \bar{a}) \bar{V} + T \mu \tilde{r}_z\\
            &&            &\leq (1- T \bar{a}) \bar{V} + T \mu \frac{\bar{a}\bar{V}}{\mu}\\
            &&            & \leq \bar{V}.
        \end{aligned}
    \end{equation}
    This implies that $x_+$ is in $X_{\bar{V}}$. Similarly, 
    due to the fact that that $E(x,z) \leq \tilde{r}_z \leq r_z$ and $x \in X_{\bar{V}}$,
    we can apply the result from Proposition
    \ref{pro:state_evolution_bound}, which shows that 
    \begin{equation}\label{equ:lem_contr_eq_1}
    \| z_+ - \bar{z}(x_+)\|  \leq \, \kappa\| z - \bar{z}(x)\| + T \gamma \| x \|
    \end{equation}
    and
    \begin{equation}
    \| z_+ - \bar{z}(x_+)\|  \leq r_z
    \end{equation}
    must hold. 
    Using Assumption \ref{ass:lyapunov_stability} in Equation \eqref{equ:lem_contr_eq_1},
    we obtain
    \begin{equation}
        \| z_+ - \bar{z}(x_+)\|  \leq \, \kappa\| z - \bar{z}(x)\| + T \hat{\gamma} \left(V(x)\right)^{\frac{1}{q}}.
    \end{equation}
    Moreover, due to \eqref{equ:ts_tilde_r_z}, we have that $\| z_+ - \bar{z}(x_+) \| \leq \tilde{r}_z$ for any 
    $T \leq T_4$. \hfill$\qed$
    \end{pf}
\end{lem}
    Lemma \ref{lem:invariance} provides invariance of $\Sigma$ for the system-optimizer dynamics 
    as well as a compact description of the interaction between the two underlying subsystems. 
    In principle, we could study the behavior of the coupled contraction estimate 
    by looking at the ``worst-case'' dynamics associated with \eqref{equ:coupled_contraction}: 
    \begin{equation}\label{equ:aux_dyn_q}
        \begin{aligned}
            && v_+ & = (1 - T \bar{a}) v + T \mu e, \\
            && e_+ & = T\hat{\gamma} v^{\frac{1}{q}}  + \kappa  e. \\
        \end{aligned}
    \end{equation}
    However, these dynamics are not Lipschitz continuous at the origin $(v,e)
    = (0,0)$ for $q>1$ and are for this reason non-trivial to analyze. 
    \par
    Nonetheless, 
    we can reformulate \eqref{equ:aux_dyn_q} such that standard tools can be used 
    to obtain important information about the behavior of the
    $V_+(T;x,z)$ and $E_+(T;x,z)$ under the combined contraction 
    established in Lemma \ref{lem:invariance} and ultimately establish asymptotic stability 
    of the system-optimizer dynamics.
\section{Asymptotic stability of the system-optimizer dynamics}\label{sec:stability}
In the following, we derive the main asymptotic stability result, which relies on 
a reformulation of the worst-case dynamics \eqref{equ:aux_dyn_q}.
\begin{prop}\label{pro:inexact_lyapunov_contraction_q}
    Let Assumptions \ref{ass:lyapunov_stability}, \ref{ass:qlip},
    \ref{ass:policy_solver}, \ref{ass:Lip_cont_z}, \ref{ass:contraction}, 
    and \ref{ass:state_evolution_bound} hold. 
    Moreover, let $\hat{\mu} \vcentcolon = L_{\phi,u} e^{T_1 L_{\phi,x}} \tilde{\mu}$. Then, for any $\xi \in \Sigma$ and
    any $T \leq T_4$, we have $\xi_+ \in \Sigma$ and the following holds: 
    \begin{equation}\label{equ:linear_coupled_ineq}
    \begin{aligned}
        &&V_+(T;x,z)^{\frac{1}{q}} &\leq (1 - T \bar{a})^{\frac{1}{q}} V(x)^{\frac{1}{q}} + T\hat{\mu} E(x,z), \\ 
        && E_+(T;x,z) & \leq  T \hat{\gamma} V(x)^{\frac{1}{q}}  + \kappa  E(x,z).
    \end{aligned}
    \end{equation}
    \begin{pf}
        The fact that $\xi_+ \in \Sigma$ is a direct consequence of Lemma \ref{lem:invariance}. 
        Moreover, due to Assumption \ref{ass:qlip}, the following holds, for any $x \in X_{\bar{V}}$: 
        \begin{equation*}
        \begin{aligned}
            &&V(\psi(T;\, &x, M_{u,z}z))^{\frac{1}{q}} \leq V(\psi(T; x, M_{u,z}\bar{z}(x))^{\frac{1}{q}} \\
            && &+ \tilde{\mu} \| \psi(T;x,M_{u,z}z) - \psi(T;x,M_{u,z}\bar{z}(x))\|
        \end{aligned}
        \end{equation*}
        and, using the nominal Lyapunov contraction and Proposition
        \ref{pro:state_evolution_bound}, we obtain, 
        for any $\xi \in \Sigma$, that
        \begin{equation*}
            V(\psi(T; x, M_{u,z}z))^{\frac{1}{q}} \leq (1-T\bar{a})^{\frac{1}{q}}V(x)^{\frac{1}{q}} + T \hat{\mu} \| z - \bar{z}(x)\|,
        \end{equation*}
        where $\hat{\mu} \vcentcolon = L_{\phi,u} e^{T_1 L_{\phi,x}} \tilde{\mu}$. \hfill$\qed$
    \end{pf}
\end{prop}
Unlike \eqref{equ:aux_dyn_q}, the worst-case dynamics associated with \eqref{equ:linear_coupled_ineq} are not only
Lipschitz continuous, but can also be cast as a linear positive system. We define the following dynamical system based on
\eqref{equ:linear_coupled_ineq}.
\begin{defn}[Auxiliary dynamics]\label{def:aux_dyn}
    We will refer to the linear time-invariant discrete-time dynamical system
    \begin{equation}\label{equ:aux_dyn}
    \begin{aligned}
        && \nu_+ & = (1 - T \bar{a})^{\frac{1}{q}} \nu + T \hat{\mu} \epsilon, \\
        && \epsilon_+ & = T \hat{\gamma} \nu  + \kappa \epsilon, \\
    \end{aligned}
    \end{equation}
    with states $\nu, \, \epsilon \in \mathbb{R}$,  as \textit{auxiliary dynamics}. 
    Due to the definitions of $\kappa, \hat{\mu}$, and $\hat{\gamma}$ and Assumption 
    \ref{ass:lyapunov_stability}, \eqref{equ:aux_dyn} is a positive system \cite{Kaczorek2008}.
\end{defn}
\begin{rem}\label{rem:upper_bound}
    Given the definition of the auxiliary dynamics in Definition \ref{def:aux_dyn}, for 
    any $\xi \in \Sigma$, if $V(x)^{\frac{1}{q}} = \nu$ and $E(x,z)=\epsilon$, then 
    $V_+(T;x,z)^{\frac{1}{q}} \leq \nu_+$ and $E_+(T;x,z) \leq \epsilon_+$. For this reason, intuitively, we can study 
    the stability of the auxiliary dynamics and infer stability properties of the 
    combined system-optimizer dynamics. This concept will be later formalized with the explicit construction of 
    a Lyapunov function for the system-optimizer dynamics in Theorem \ref{thm:main_thm}.
\end{rem}
We exploit properties of 
positive systems in order to construct an explicit linear Lyapunov function for 
the auxiliary dynamics which can be rewritten in the compact form
\begin{equation}\label{equ:possys}
    w_+ = A_a w,
\end{equation}
where 
\begin{equation}
    A_a \vcentcolon=
    \begin{bmatrix}
        (1 - T \bar{a})^{\frac{1}{q}} &\quad T \hat{\mu} \\ 
        T \hat{\gamma}   &\quad \kappa
    \end{bmatrix},
\end{equation}
and $w\vcentcolon = (\nu, \epsilon)$. We will make use of the following result adapted from \cite{Kaczorek2008}. 
\begin{thm}[Stability of positive systems]\label{thm:ps_stab}
    A positive discrete-time linear system of the form 
    \begin{equation}\label{equ:ps}
        w_+ = A w,
    \end{equation}
    where $A \in \mathbb{R}_{\geq 0}^{n \times n}$ and $w \in \mathbb{R}_{\geq 0}^n$, is asymptotically stable 
    if there exist a strictly positive vector $\hat{w} \in \mathbb{R}^{n}_{>0}$ 
    and a strictly positive constant $\hat{d}$ such that 
    \begin{equation}
        \underset{i=1,\dots,n}{\max} \, [\, (A^{\top} - \mathbb{I})\hat{w} \,]_i \leq - \,\hat{d}.
    \end{equation}
    Moreover, the linear function $V_l(w) \vcentcolon = \hat{w}^\top w$ is a 
    Lyapunov function for \eqref{equ:ps} in $\mathbb{R}_{\geq0}^n$ and the following holds:
    \begin{equation}
        V_l(w_+) - V_l(w) \leq -\hat{d} \cdot \| w \|. 
    \end{equation}

    \begin{pf}
        See Appendix \ref{app:ps_stab}.\hfill$\qed$
    \end{pf}
\end{thm}
\begin{figure}
\centering
\includegraphics[scale=1.0]{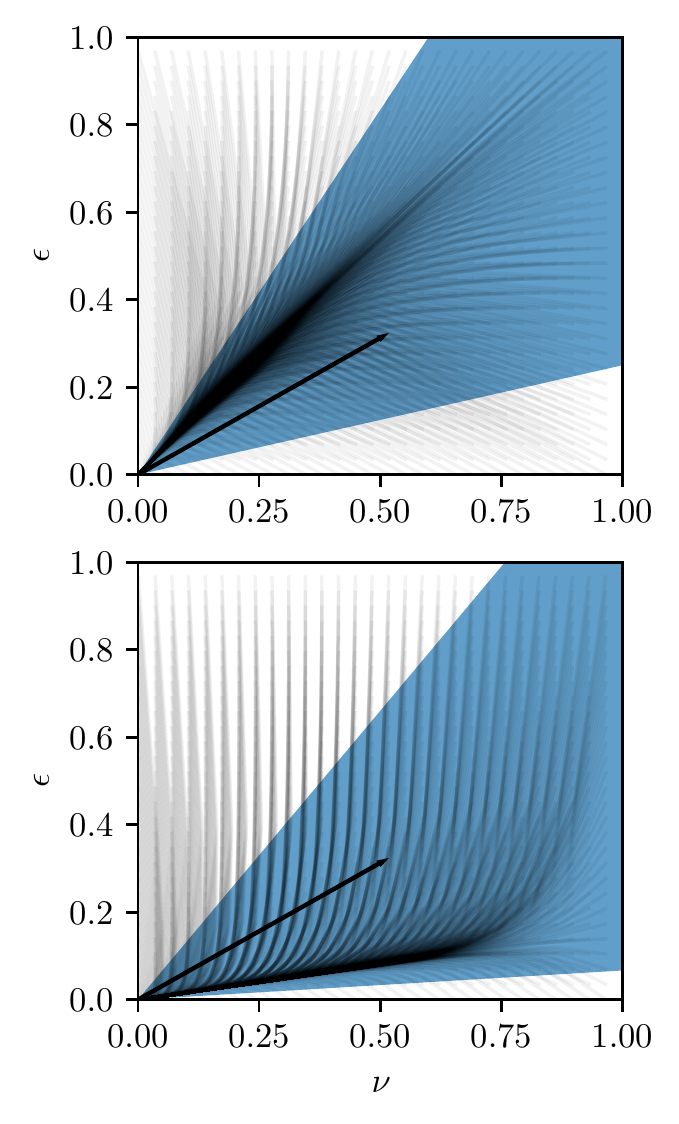}
\caption{Trajectories of the auxiliary dynamics \eqref{equ:aux_dyn} for different 
    initial conditions ($\kappa=0.4$, $\bar{a}=0.5$, $\hat{\gamma}=0.2$, 
    $\hat{\mu}=0.1$) - $T = 1.0$ (top) and $T = 0.4$ (bottom). The black 
    vector describes the direction defined by $\hat{w}$ as in Theorem 
    \ref{thm:stability_aux}, while the shaded area defines the cone that 
    contains all the vectors that would satisfy \eqref{equ:ineq_linear_lyapunov}, 
    i.e., all the vectors $\hat{w}$ that define a valid Lyapunov function $V_l(w) = \hat{w}^{\top}w$.} 
    \label{fig:trajectories}
\end{figure}
\begin{thm}\label{thm:stability_aux}
    The positive discrete-time linear system \eqref{equ:possys} is asymptotically stable if 
    and only if the following condition is satisfied:
    \begin{equation}\label{equ:aux_stability_ts_condition}
        T^2 \hat{\mu}\hat{\gamma}- (1-\kappa) (1 - (1 - T \bar{a})^{\frac{1}{q}}) < 0,
    \end{equation}
    which holds for any sufficiently small sampling time $T \leq T_5 \vcentcolon = \frac{\beta (1 - \kappa)}{\hat{\mu}}$.
    Moreover, the function $V_l(w) \vcentcolon= \hat{w}^{\top}w$, where 
    \begin{equation}
        \hat{w} = \begin{bmatrix} 1 \\ \beta\end{bmatrix}, \quad \MR{with} \quad \beta \vcentcolon = \frac{1}{2}\frac{\bar{a}}{q\hat{\gamma}},
    \end{equation}
    is a Lyapunov function for \eqref{equ:aux_dyn} in
    $\mathbb{R}_{\geq0}^2$.
    \begin{pf}
        In order to prove that $\hat{w}^{\top}w$ is a Lyapunov function for \eqref{equ:aux_dyn}
        it suffices to apply Theorem \ref{thm:ps_stab}. The system of inequalities
        \begin{equation}\label{equ:ineq_linear_lyapunov}
            \begin{aligned}
                &&-1 + (1-T \bar{a})^{\frac{1}{q}} \,  + T \hat{\gamma} \,\beta &< 0, \\
                &&T \hat{\mu}\,  + (\kappa - 1) \, \beta &< 0, \\
                &&\beta &> 0
            \end{aligned}
        \end{equation}
    admits a solution if 
    \begin{equation}\label{equ:ps_con}
        T \frac{\hat{\mu}}{1-\kappa} < \beta <  \frac{1 - (1 - T \bar{a})^{\frac{1}{q}}}{T \hat{\gamma}}.
    \end{equation}
    This condition can always be satisfied for a sufficiently small $T$. In fact, 
    it is easy to show that the limits for $T\to 0$ of the upper and lower bounds on 
    $\beta$ are 0 and $\frac{\bar{a}}{q\hat{\gamma}} > 0$, respectively, such that, by continuity, 
    there must exist a strictly positive constant $T_5$ such that \eqref{equ:aux_stability_ts_condition} 
    is satisfied for any $T \leq T_5$.
    However, in order to make $\beta$ independent of $T$, we note
    that, due to convexity, $1 - (1 - T \bar{a})^{\frac{1}{q}} \geq \frac{\bar{a}}{q}T$ 
    for any $T \geq 0$. Using this lower bound we can simplify the upper bound in \eqref{equ:ps_con} as
    \begin{equation}
        \beta < \frac{\frac{\bar{a}}{q}T}{\hat{\gamma}T} = 
        \frac{\bar{a}}{q\hat{\gamma}} \leq \frac{1 - (1 - T \bar{a})^{\frac{1}{q}}}{T \hat{\gamma}}.
    \end{equation}
    Choosing $\beta$ to be half of this upper bound, i.e., $\beta \vcentcolon = \frac{1}{2}\frac{\bar{a}}{q\hat{\gamma}}$, we obtain that \eqref{equ:ps_con} is 
    satisfied for any $T \leq T_5 \vcentcolon = \frac{\beta (1 - \kappa)}{\hat{\mu}}$, which concludes the proof. \hfill$\qed$
    \end{pf}
\end{thm}
Theorem \ref{thm:stability_aux} shows that (exponential) asymptotic stability of the
auxiliary dynamics holds under the condition that the sampling time $T$ satisfies 
inequality \eqref{equ:aux_stability_ts_condition} for given $\hat{\mu}$, $\bar{a}$, $\hat{\gamma}$, and $\kappa$. 
Figure \ref{fig:trajectories} illustrates the meaning of Theorem \ref{thm:stability_aux} 
by showing the trajectories of the auxiliary system in a phase plot for 
fixed values of the parameters $\hat{\mu}, \bar{a}, \kappa$, and $\hat{\gamma}$, two different values 
of the sampling time $T$ and for different initial conditions.
In the following, we establish the main result of the section by exploiting the Lyapunov decrease 
established in Theorem \ref{thm:stability_aux} for the auxiliary dynamics to construct a Lyapunov function 
for the combined system-optimizer dynamics \eqref{equ:so_dyn_compact}.
\andrea{remark on $\hat{\kappa}$ and number of iterations?!}
\begin{thm}\label{thm:main_thm}
    Let Assumptions \ref{ass:lyapunov_stability},
    \ref{ass:qlip}, \ref{ass:policy_solver}, \ref{ass:Lip_cont_z},
    \ref{ass:contraction}, and \ref{ass:state_evolution_bound} hold. Then, for any $T \leq \min \{ T_4, T_5\}$, 
    the origin is an exponentially asymptotically stable equilibrium with region of attraction 
    $\Sigma$ for the combined system-optimizer dynamics \eqref{equ:so_dyn_compact}. 
    In particular, the function
    \begin{equation}
        V_{\mathrm{so}}(\xi) \vcentcolon = \hat{w}^{\top}\begin{bmatrix} V(x)^{\frac{1}{q}} \\ \| z - \bar{z}(x) \|\end{bmatrix},
    \end{equation}
    where $\hat{w}$ is defined according to Theorem \ref{thm:stability_aux},
    is a Lyapunov function in $\Sigma$ for the system \eqref{equ:so_dyn_compact} and the origin $(x,z) = \xi = 0$. 
    \begin{pf}
        We can derive an upper bound for $V_{\mathrm{so}}(\xi)$ as follows:
        \begin{equation}
            \begin{aligned}
            &&V_{\mathrm{so}}(\xi) &=  V(x)^{\frac{1}{q}} + \beta \| z - \bar{z}(x) \| \\
            && &\leq  a_2^{\frac{1}{q}} \| x \| + \beta \| z - \bar{z}(x) \| \\
            && &\leq  a_2^{\frac{1}{q}} \| x \| + \beta (\| z \| + \| \bar{z}(x) \| )\\
            && &\leq ( a_2^{\frac{1}{q}} + \sigma \beta)\| x \| + \beta \| z \| \\
            && &\leq \underbrace{\max\{ a_2^{\frac{1}{q}} + \sigma \beta, \beta\}}_{
            =\vcentcolon\tilde{w}_2} \cdot \left(\| x \| + \| z \|\right) \\
            && &\leq \tilde{w}_2 \cdot (\| x\|_1 + \| z \|_1) = \tilde{w}_2 \cdot 
                \| \xi \|_1 \\
            && &\leq \tilde{w}_2 \sqrt{n_x + n_z} \cdot \| \xi \|. \\
            \end{aligned}
        \end{equation}
        In order to derive a lower bound, we proceed as follows. 
        Using the reverse triangle inequality and Lipschitz continuity 
        of $\bar{z}(x)$, we obtain
        \begin{equation}
        \begin{aligned}
            &&V_{\mathrm{so}}(\xi) &\geq  V(x)^{\frac{1}{q}}+ \beta (\| z \| - \sigma \| x \|) \\ 
            && &\geq  a_1^{\frac{1}{q}} \| x \|+ \beta (\| z \| - \sigma \| x \|) \\ 
            &&              &= ( a_1^{\frac{1}{q}} - \beta \sigma ) \| x \| + \beta \| z \|. 
        \end{aligned}
        \end{equation}
        If $ a_1^{\frac{1}{q}} - \beta \sigma > 0$, then we can readily compute a lower bound:
        \begin{equation}\label{equ:main_thm_bound_1}
        \begin{aligned}
            &&V_{\mathrm{so}}(\xi) &\geq ( a_1^{\frac{1}{q}} - \beta \sigma) \| x \| 
            + \beta \| z \|  \\
            && &\geq \frac{ a_1^{\frac{1}{q}} - \beta \sigma}{\sqrt{n_x}} \| x \|_1 
            + \frac{\beta}{\sqrt{n_z}}
            \| z \|_1  \\
            && &\geq \underbrace{\min \bigg\{ \frac{ a_1^{\frac{1}{q}} 
            - \beta \sigma}{\sqrt{n_x}}, \frac{\beta}{\sqrt{n_z}}
    \bigg\}}_{=\vcentcolon\tilde{w}_{1,1}} (\| x \|_1 
            + \| z \|_1 ) \\
            && & \geq \tilde{w}_{1,1} \cdot \| \xi \|.
        \end{aligned}
        \end{equation}
        If instead $ a_1^{\frac{1}{q}} - \beta \sigma \leq 0$, we define the auxiliary lower bound
        % \begin{equation}
        %     \begin{aligned}
        %     \hat{V}_{\mathrm{so}}(R,z) \vcentcolon=
        %     \,\,\,&\underset{x}{\min}	&& ( a_1^{\frac{1}{q}} - \sigma \beta) \| x \| + \beta \| z \| \\
        %     &\,\, \MR{s.t.}  &&\| x \| = R.\\
        %     \end{aligned}
        % \end{equation}
        \begin{equation}
            \hat{V}_{\mathrm{so}}(x,z) \vcentcolon=  a_1^{\frac{1}{q}} \| x\| + \beta \| z - \bar{z}(x)\|. 
        \end{equation}
        Since $V_{\mathrm{so}}(\xi) \geq \hat{V}_{\mathrm{so}}(x,z)$, if we can show that 
        $\hat{V}_{\mathrm{so}}(x,z)$ can be lower bounded by a properly constructed function of $z$, we can use this function to 
        lower bound $V_{\mathrm{so}}(\xi)$ too. To this end, we first observe that, since we assumed that 
        $ a_1^{\frac{1}{q}} - \beta \sigma \leq 0$, for any $x$ such that $\| x \| \leq \frac{1}{\sigma} \| z\|$, we 
        have that 
        \begin{equation}
            \begin{aligned}
            &&\hat{V}_{\mathrm{so}}(x,z) &\geq ( a_1^{\frac{1}{q}} - \beta \sigma ) \| x \| + \beta \| z \| \\
            && &\geq \underset{x \,\MR{s.t.} \,\| x \| \leq \frac{1}{\sigma} \| z \|}{\min}( a_1^{\frac{1}{q}} - \beta \sigma ) \| x \| + \beta \| z \| \\
            && &\geq \frac{a_1^{\frac{1}{q}} }{\sigma} \| z \|,
            \end{aligned}
        \end{equation}
        where, for the minimization, we have used the fact that the objective $( a_1^{\frac{1}{q}} - \beta \sigma ) \| x \| + 
        \beta \| z\|$ is monotonically nonincreasing in $\| x \|$ such that the minimum is attained at 
        the boundary of the interval for $\| x \| = \frac{1}{\sigma} \| z\|$.
        Similarly, 
        for any $x$ such that $\| x\| \geq \frac{1}{\sigma} \| z \|$, we can use the fact that 
        \begin{equation*}
            \begin{aligned}
                &&\hat{V}_{\mathrm{so}}(x,z) &=  a_1^{\frac{1}{q}} \| x\| + 
                \beta \| z - \bar{z}(x)\| 
                \geq  a_1^{\frac{1}{q}} \| x \| \geq \frac{ a_1^{\frac{1}{q}} }{\sigma} \| z \|.
            \end{aligned}
        \end{equation*}
        Hence, we can conclude that 
        \begin{equation}
            V_{\mathrm{so}}(\xi) \geq \hat{V}_{\mathrm{so}}(x,z)  \geq \frac{ a_1^{\frac{1}{q}} }{\sigma} \| z \|
        \end{equation}
        for any $x$.
        Summing this last inequality and $V_{\mathrm{so}}(\xi)  \geq  a_1^{\frac{1}{q}} \| x \|$, 
        we obtain
        \begin{equation}
        \begin{aligned}
            &&V_{\mathrm{so}}(\xi) &\geq \frac{ a_1^{\frac{1}{q}} }{2} \| x \|
            + \frac{ a_1^{\frac{1}{q}} }{2\sigma} \| z \|  \\
            && &\geq \frac{ a_1^{\frac{1}{q}}}{2 \sqrt{n_x}} \| x \|_1 
            + \frac{ a_1^{\frac{1}{q}}}{2\sigma\sqrt{n_z}}
            \| z \|_1  \\
            && &\geq \underbrace{\min \bigg\{ \frac{ a_1^{\frac{1}{q}} 
            }{2\sqrt{n_x}}, \frac{ a_1^{\frac{1}{q}}}{2\sigma\sqrt{n_z}}\bigg\}}_{
        =\vcentcolon\tilde{w}_{1,2}} (\| x \|_1 
            + \| z \|_1 ) \\
            && & \geq \tilde{w}_{1,2} \cdot \| \xi \|.
        \end{aligned}
        \end{equation}
        Together with \eqref{equ:main_thm_bound_1}, we can define 
        \begin{equation}
            \tilde{w}_1 \vcentcolon = 
            \begin{cases}
                \tilde{w}_{1,1}, \quad \text{if} \quad a_1^{\frac{1}{q}} - \beta \sigma > 0, \\
                \tilde{w}_{1,2}, \quad \text{otherwise}, \\
            \end{cases}
        \end{equation}
        and conclude that $V_{\MR{so}}(\xi) \geq \tilde{w}_1 \cdot \| \xi\|$. 
        Finally, the Lyapunov decrease can be derived. For
        given $x$ and $z$, let $\epsilon = E(x,z)$ and $\nu =
        V(x)^{\frac{1}{q}}$. Then the following holds.
        \begin{equation}
            \begin{aligned}
                &&V_{\mathrm{so}}(\xi_+) =&\,V_+(T;x,z)^{\frac{1}{q}} +
                \beta E_+(T;x,z) \\
                && \stackrel{\mathclap{\rm Remark
        \,\,\ref{rem:upper_bound}}}{\leq}&\,
                \nu_+ + \beta \epsilon_+ \\
                && \stackrel{\mathclap{\rm Theorem \,\,
        \ref{thm:ps_stab}}}{\leq}&\,
                    \nu + \beta \epsilon - \hat{d} \cdot \| 
                    (\nu, \epsilon)\|_1\\
                    && =&\,  V(x)^{\frac{1}{q}} + \beta E(x,z) - \hat{d} \cdot
                    \| (\nu, \epsilon) \|_1 \\
                    && =&\, V_{\mathrm{so}}(\xi) - \hat{d}
                    \cdot (V(x)^{\frac{1}{q}} + \| z - \bar{z}(x)\|),  \\ 
            \end{aligned}
        \end{equation}
        where we have used the $\ell_1$ norm due to the intermediate result in the proof 
        of Theorem \ref{thm:ps_stab}.
        Let $\Delta V_{\mathrm{so}} (\xi)\vcentcolon=-\hat{d} \cdot (
        V(x)^{\frac{1}{q}} +  \| z - \bar{z}(x)\|)$
        denote the Lyapunov decrease.
        Using the same procedure used to derive the lower bound for $V_{\mathrm{so}}(\xi)$, 
        we can show that, if $a_1^{\frac{1}{q}} - \sigma > 0$, we can write 
        \begin{equation*}
        \begin{aligned}
            &&\Delta V_{\mathrm{so}}(\xi) &\leq -\hat{d} \cdot \left((a_1^{\frac{1}{q}} - \sigma) \| x \| 
            + \| z \|\right)  \\
            && &\leq -\underbrace{\hat{d} \cdot \min \bigg\{ \frac{ a_1^{\frac{1}{q}} - \sigma
}{\sqrt{n_x}}, \frac{1}{\sqrt{n_z}} \bigg\}}_{=\vcentcolon\tilde{w}_{3,1}} (\| x \|_1 
            + \| z \|_1 ) \\
            && & \leq -\tilde{w}_{3,1} \cdot \| \xi \|.
        \end{aligned}
        \end{equation*}
        Else, if $a_1^{\frac{1}{q}} - \sigma \leq 0$, we can obtain the following bound:
        \begin{equation*}
            \begin{aligned}
                &&\Delta V_{\mathrm{so}} (\xi)&=-\hat{d}
                \cdot \left(V(x)^{\frac{1}{q}} +  \| z -
                \bar{z}(x)\|\right) 
                \leq -\frac{\hat{d} a_1^{\frac{1}{q}}}{\sigma} \| z \|.
        \end{aligned}
        \end{equation*}
        Summing this last inequality and $\Delta V_{\mathrm{so}} (\xi) \leq -\hat{d} a_1^{\frac{1}{q}} \| x \|$, 
        we obtain
        \begin{equation*}
        \begin{aligned}
            &&\Delta V_{\mathrm{so}}(\xi) &\leq 
            -\frac{\hat{d}a_1^{\frac{1}{q}} }{2} \left(\| x \| +
            \frac{1}{\sigma} \| z \|\right) \\ 
            && & \leq -\frac{\hat{d}a_1^{\frac{1}{q}} }{2}\left(\frac{1}{
                \sqrt{n_x}} \| x \|_1 
            + \frac{1}{\sigma\sqrt{n_z}}\| z \|_1\right)
              \\
            && &\leq -\underbrace{\frac{\hat{d}a_1^{\frac{1}{q}}}{2}
              \cdot
              \min \bigg\{ \frac{1}{\sqrt{n_x}},
              \frac{1}{\sigma\sqrt{n_z}}\bigg\}}_{
          =\vcentcolon\tilde{w}_{3,2}} (\| x \|_1 
            + \| z \|_1 ) \\
            && & \leq -\tilde{w}_{3,2} \cdot \| \xi \|.
        \end{aligned}
        \end{equation*}
        We define 
        \begin{equation}
            \tilde{w}_3 \vcentcolon = 
            \begin{cases}
                \tilde{w}_{3,1}, \quad \text{if} \quad a_1^{\frac{1}{q}} - \sigma > 0, \\
                \tilde{w}_{3,2}, \quad \text{otherwise}, \\
            \end{cases}
        \end{equation}
        and conclude that $\Delta V_{\MR{so}}(\xi) \leq -\tilde{w}_3 \cdot \| \xi\|$. 
        Hence, we can define the $\mathcal{K}_{\infty}$ functions 
        $\alpha_{\mathrm{so},1}(\| \xi \|)\vcentcolon= \tilde{w}_1 \cdot \| \xi \|$ 
        and $\alpha_{\mathrm{so},2}(\| \xi \|)\vcentcolon=\tilde{w}_2 \cdot \| \xi \|$
        and the positive definite function $\alpha_{\mathrm{so},3}(\| \xi \|)\vcentcolon=\tilde{w}_3 \cdot \| \xi \|$,
        such that
        \begin{equation}
        \begin{aligned}
            \alpha_{\mathrm{so},1}(\| \xi \|) \leq V_{\mathrm{so}}(\xi) \leq \alpha_{\mathrm{so},2}(\| \xi \|) \\
            V_{\mathrm{so}}(\xi_+) - V_{\mathrm{so}}(\xi) \leq -\alpha_{\mathrm{so},3}(\| \xi \|),
        \end{aligned}
        \end{equation}
        i.e., $V_{\mathrm{so}}(\xi)$ is a Lyapunov function in $\Sigma$ for the system-optimizer dynamics $\xi_+ = \Phi(T;\xi)$
        and the equilibrium $\xi=0$, for any $T \leq \min\{T_4, T_5\}$. \hfill$\qed$
    \end{pf}
\end{thm}
Theorem \ref{thm:main_thm} shows that, for a sufficiently short sampling time $T$, the system-optimizer dynamics in 
\eqref{equ:so_dyn} are asymptotically exponentially stable and that $V_{\text{so}}$ is a Lyapunov function in $\Sigma$.
We observe that the obtained Lyapunov function is a positive linear combination of the ideal NMPC Lyapunov function $V$ 
and the error $\| z - \bar{z}(x)\|$ (which in turn is a Lyapunov function for the error dynamics). Loosely speaking, depending on the value 
of $\beta$, $V_{\text{so}}$ gives more ``importance'' to either $V$ (for small values of $\beta$) or $\| z - \bar{z}(x)\|$ (for large values of $\beta$).
\par
In the next section, in order to illustrate Theorem \ref{thm:main_thm}, we 
use a variant of the classical example from \cite{Chen1998}. 
\section{Illustrative example}
\begin{figure}[t]
\centering
\includegraphics[scale=0.9]{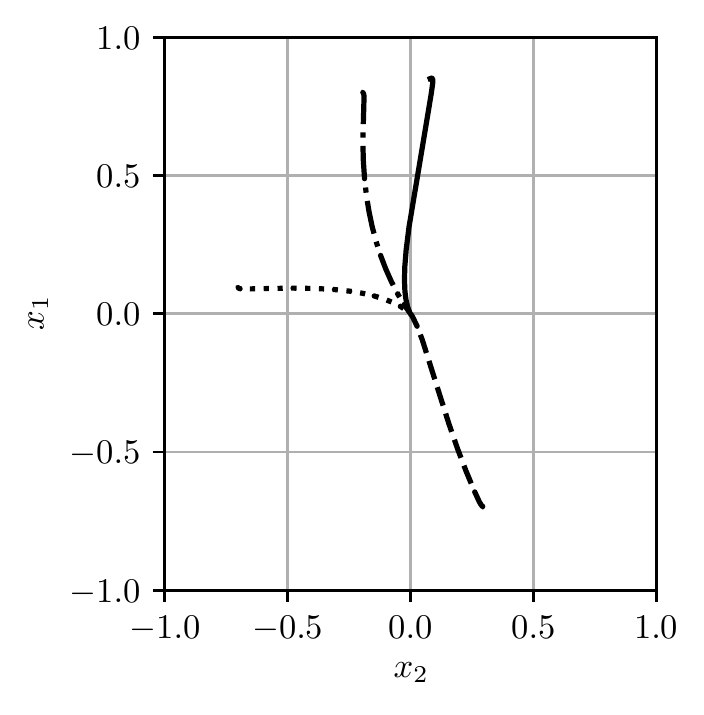}
\caption{Illustrative example adapted from \cite{Chen1998} - closed-loop state trajectories 
obtained using the approximate feedback policy computed with a single iteration of a Gauss-Newton real-time algorithm.
} 
\label{fig:chen_states}
\end{figure}
\begin{figure}[t]
\centering
\includegraphics[scale=1.0]{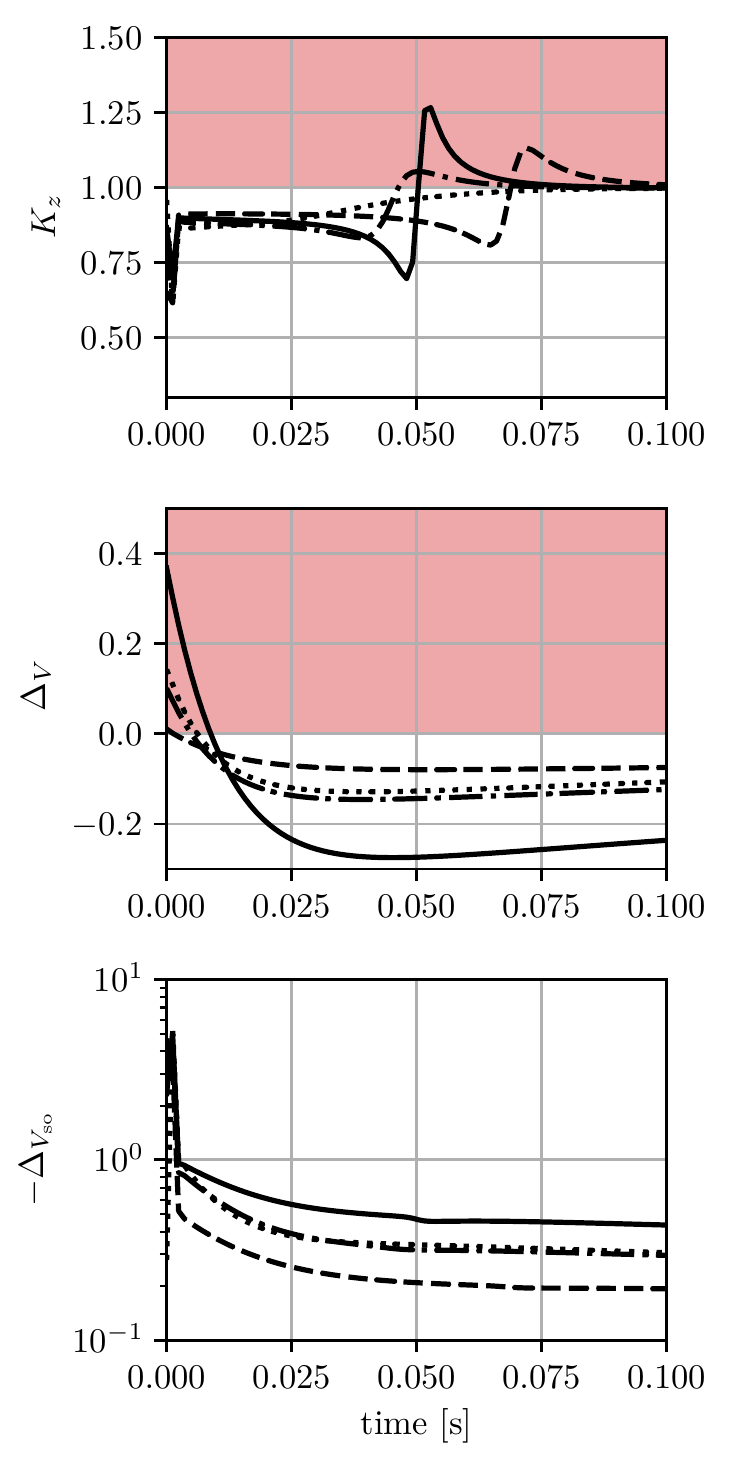}
\caption{Illustrative example adapted from \cite{Chen1998} - although 
the numerical error and the value function do not necessarily decrease monotonically over time, the Lyapunov function for 
the combined system-optimizer dynamics $V_{\mathrm{so}}(\xi)$ does decrease.} 
\label{fig:chen_lyapunov}
\end{figure}
We regard an optimal control problem of the form in \eqref{eq:nmpc} and define the continuous-time dynamics as
\begin{equation}\label{equ:chen_dyn}
    \phi(x,u) \vcentcolon= \begin{bmatrix} x_2 + u(\mu + (1-\mu)x_2) \\ 
                                          x_1 + u(\mu - 4(1-\mu)x_2)\end{bmatrix}.
\end{equation}
In order to compute an LQR-based terminal cost, the linearized dynamics are defined as
\begin{equation}
    A_{\MR{c}} \vcentcolon = \frac{\partial \phi}{\partial x}(0,0), \quad B_{\MR{c}} \vcentcolon = \frac{\partial \phi}{\partial u}(0,0), 
\end{equation}
and discretized using exact discretization with piece-wise constant 
parametrization of the control trajectories:
\begin{equation*}
    A \vcentcolon= \exp{(A_{\MR{c}} T_{\MR{d}})}, \quad B \vcentcolon= \left(\int_{\tau = 0}^{T_{\MR{d}}}\exp{(A_{\MR{c}}\tau) \MR{d}\tau}\right)B_{\MR{c}},
\end{equation*}
where $T_{\MR d}$ denotes the discretization time.
We compute the solution $P$ to the discrete-time 
algebraic Riccati equation 
\begin{equation*}
        P = A^{\top}PA - (A^{\top}PB)(R + B^{\top}PB)^{-1}(B^{\top}PA) + Q,
\end{equation*}
where $Q = 0.1 \cdot \mathbb{I}_2$ and $R = 0.1$ such that the cost functionals can be defined as 
\begin{equation}
    L(x,u) \vcentcolon = x^{\top}Qx + u^{\top}Ru, \quad m(x) \vcentcolon = x^{\top}Px,
\end{equation}
and we impose simple bounds on the input $-2 \leq u \leq 2$.
\par
We set the prediction horizon $T_{\MR{f}}=0.3$ and discretize the resulting continuous-time optimal control problem using 
direct multiple shooting with $N=5$ shooting nodes. The Euler discretization is used for the 
cost integral and explicit RK4 is used to discretize the dynamics. In order to 
solve the resulting discretized optimal control problem, 
we use the standard RTI approach, with Gauss-Newton iterations and a fixed Levenberg-Marquardt-type term. 
A single SQP iteration per sampling time is carried out.
\par
In order to compute an estimate for the constants involved in the definition of the 
Lyapunov function in Theorem \ref{thm:main_thm}, we regard six different initial conditions, 
and control the system using the feedback policy associated with the exact solution to the discretized optimal control problem. 
For every state $x$ in the obtained trajectories, we evaluate the optimal cost $V(x)$ and the primal-dual optimal solution 
$\bar{z}(x)$. With these values, we can estimate the constants $a_1, a_2, a_3$ in Assumption \ref{ass:lyapunov_stability}, 
the constant $\tilde{\mu}$ in Assumption \ref{ass:qlip} and the constant $\sigma$ in Theorem \ref{lem:contraction}. 
Moreover, for any state visited, we carry out a limited number of iterations of the optimizer in 
order to estimate the contraction rate $\hat{\kappa}$. Choosing the sampling time $T=0.0012$, we obtain 
the estimates $\kappa = 0.882, \bar{a} = 1.157, \hat{\gamma} = 70.23$, and $\hat{\mu} = 0.360$, such that
the parameter involved in the definition of the Lyapunov function for the combined system-optimizer dynamics takes the value $\beta = 0.0041$ and 
we have $T_5 = 0.037 \geq T$.
All the computations were carried out using CasADi \cite{Andersson2018} and its interface to Ipopt \cite{Waechter2009} and the code for 
the illustrative example is made available at \url{https://github.com/zanellia/nmpc\_system_optimizer\_lyapunov}.
\par
Figure \ref{fig:chen_states} shows the state trajectories obtained controlling the system 
using the approximate feedback law starting from the selected initial conditions. 
Figure \ref{fig:chen_lyapunov} shows the behavior of $\| z - \bar{z}(x)\|$, $V(x)$ and 
 $V_{\mathrm{so}}(\xi)$ over time through the compact metrics 
 \begin{equation}
     K_z \vcentcolon = \frac{\| z_+ - \bar{z}(x_+)\|}{\|z - \bar{z}(x)\|},
 \end{equation}
 \begin{equation}
     \Delta_{V} \vcentcolon = \frac{V(x_+) - V(x)}{T_s},
 \end{equation}
 \begin{equation}
 \Delta_{V_{\text{so}}} \vcentcolon = \frac{V_{\text{so}}(\xi_+) - V_{\text{so}}(\xi)}{T_s}.
 \end{equation}
In particular, Figure \ref{fig:chen_lyapunov} shows that, although the numerical error $\| z - \bar{z}(x)\|$ and the 
value function does not necessarily decrease over time, the constructed Lyapunov function $V_{\mathrm{so}}(\xi)$ does decrease.

\section{Conclusions}
In this paper, we presented a novel asymptotic stability results for inexact MPC  
relying on a limited number of iterations of an optimization algorithm.
A class of optimization methods with Q-linearly convergent iterates has been regarded and, under the assumption that 
the ideal feedback law is stabilizing, we constructed a Lyapunov function for the system-optimizer dynamics. 
These results extend the attractivity proofs present in the literature 
which rely instead on the assumption that inequality constraints are either absent or inactive in the attraction region considered.
Moreover, with respect to more general results on suboptimal MPC (cf. \cite{Scokaert1999}, \cite{Pannocchia2011}, \cite{Allan2017}), 
we analyzed the coupled system and optimizer dynamics and proved its asymptotic stability.
\par
Future research will investigate how to relax the requirements of Lipschitz continuity of the solution localization $\bar{z}$ and of Q-linear convergence.
\color{black}
\begin{comment}

\section{Open Question}
We have the following coupled Lyapunov-error contractions
    \begin{equation}\label{eq:contr_Lyapunov}
    \begin{aligned}
        \!\!\!\!V_+ \leq & \,\, \left(\frac{M - m}{M}\right)V  \\
                                  & + \mu \left(\frac{\sqrt{V}}{\sqrt{m}}\norm{z - \bar{z}(x)} + \norm{z - \bar{z}(x)}^2\right)
    \end{aligned} 
    \end{equation}
and
\begin{equation}\label{eq:contr_z}
    \| z_+ - \bar{z}(x_+)\|  \leq \, \delta\| z - \bar{z}(x)\| + \frac{\gamma}{\sqrt{m}} \sqrt{V}, 
\end{equation}
where $V \vcentcolon= V(x)$.
\subsection{Proof sketch}
\begin{itemize}
    \item assume that $\| z - \bar{z}(x) \| \leq \rho \sqrt{V}$ for some $\rho > 0$. 
    \item show that, for sufficiently small $\rho$, $V^{j} \leq V$ and $\| z^j - \bar{z}^j \| \leq \| z - \bar{z} \|$, 
        for any $j \geq k$.
    \item show that $V^j$ will decrease below a certain value $\kappa_{V} V$ with $\kappa_{V} < 1$.
    \item show that for the decreased value $V^j \leq \kappa_{V} V$, we can ensure that $\| z^j - \bar{z}^j\|$ will 
        also decrease
\end{itemize}
\end{comment}
\color{black}
\bibliographystyle{plain}
\bibliography{syscop} 
\appendix
\section{Proof of Proposition 3.2.16}\label{app:g_lemma}
    In the following, we will use the shorthand $u=M_{u,z}z$. 
    First, notice that  
    \begin{equation}\label{equ:app_sol}
        \psi(t; x, u) = x +  \int_0^{t} \phi(\psi(\tau; x, u), u)\, \MR{d}\tau, \\
    \end{equation}
    and that, due to continuity of $\psi(t; x, u)$ with respect to $t$, 
    \andrea{would love to justify this more thoroughly...}
    for all $x \in X_{\bar{V}}$ and all $z \in \mathcal{B}(\bar{z}(x), r_z)$, there 
exists a $T' > 0$, such that, for all $t \leq T'$, we have $\psi(t;x,u) \in X_{\bar{V}}\oplus\mathcal{B}(0,\rho)$. 

Using \eqref{equ:app_sol} and Assumption \ref{ass:state_evolution_bound}, we obtain
    that, for any $x \in X_{\bar{V}}$, all $z \in \mathcal{B}(\bar{z}(x), r_z)$ and 
    all $T \leq T'$, the following holds:
    \begin{equation*}
    \begin{aligned}
        &&\|\psi(T; x, u) - x\|     &= \Big\|\int_0^{T}  \phi(\psi(\tau; x, u), u)\, \MR{d}\tau \Big\|\\
        &&                          &\leq \int_0^{T} \| \phi(\psi(\tau; x, u), u)\|\, \MR{d}\tau \\
        &&                          &= \int_0^{T} \| \phi(\psi(\tau; x, u), u) - \phi(0,0)\|\, \MR{d}\tau \\
    \end{aligned}
    \end{equation*}
    and, using the fact that $\phi(0,0) = 0$,
    \begin{equation*}
    \begin{aligned}
        &&\!\!\!\!\!\!\|&\psi(T; x, u) - x\| \\
        &&                          &\leq \int_0^{T} \!\!\left(L_{\phi,x} \| \psi(\tau; x, u)\|  + L_{\phi,u} \| u \|\right) \MR{d}\tau \\ 
        &&                          &= \int_0^{T} \!\!\left(L_{\phi,x} \| \psi(\tau; x, u) + x - x\| + L_{\phi,u} \| u \|\right) \MR{d}\tau \\ 
        &&                          &\leq \int_0^{T} \!\!\left(L_{\phi,x} \| \psi(\tau; x, u) - x\| + L_{\phi,u} \| u \| + L_{\phi,x} \| x \|\right)\MR{d}\tau \\ 
        &&                          &= \int_0^{T} \!\!L_{\phi,x} \| \psi(\tau; x, u) - x\| \,\MR{d}\tau \\
        && &\quad\quad\quad+ T (L_{\phi,u} \| u \| + L_{\phi,x} \| x \|). \\ 
    \end{aligned}
    \end{equation*}
    We can then apply the integral form of Gr{\"o}nwall's Lemma in order to obtain
    \begin{equation*}
    \begin{aligned}
        &&\|\psi(T; x, u) - x\| & \leq  T e^{L_{\phi,x}T}(L_{\phi,u} \| u \| + L_{\phi,x} \| x \|)
    \end{aligned}
    \end{equation*}
    for any $x \in X_{\bar{V}}$ and all $z \in \mathcal{B}(\bar{z}(x), r_z)$.
    Similarly, in order to prove the second inequality, we first notice that 
    there must exist a $T'' > 0$ such that, for all $x \in X_{\bar{V}}$, for all 
    $u = M_{u,z}z$ and all $u' = M_{u,z}z'$ such that $z,z' \in \mathcal{B}(\bar{z}(x), r_z)$ and for all $t \leq T''$, 
we have $\psi(t;x,u), \psi(t;x,u') \in X_{\bar{V}} \oplus \mathcal{B}(0,\rho)$. Hence, 
    for any $x \in X_{\bar{V}}$ and all $z,z' \in \mathcal{B}(\bar{z}(x), r_z)$ and 
    all $T \leq T''$ we can proceed as follows:
    \begin{equation*}
    \begin{aligned}
        &&\|\psi&(T; x, u') - \psi(T; x, u)\| \\
        && &\leq \int_0^{T} \| \phi(\psi(\tau; x, u'), u') - \phi(\psi(\tau; x, u), u)\|\, \MR{d}\tau \\
        && &= \int_0^{T} \| \phi(\psi(\tau; x, u'), u') - \phi(\psi(\tau; x, u), u) \\
        && & \quad+ \phi(\psi(\tau; x, u'), u) - \phi(\psi(\tau; x, u'), u) \|\, \MR{d}\tau \\
        && &\leq \int_0^{T} \| \phi(\psi(\tau; x, u'), u') - \phi(\psi(\tau; x, u'), u) \| \MR{d}\tau \\
        && & \quad + \int_0^{T}  \| \phi(\psi(\tau; x, u), u) - \phi(\psi(\tau; x, u'), u) \|\, \MR{d}\tau\\
        && & \leq T L_{\phi,u} \| u' - u \| \\
        && & \quad + \int_0^{T} L_{\phi,x} \| \psi(\tau; x, u') - \psi(\tau; x, u) \|\, \MR{d}\tau\\
    \end{aligned}
    \end{equation*}
    and, applying Gr{\"o}onwall's Lemma, we obtain
    \begin{equation*}
    \begin{aligned}
        &&\|\psi(T; x, u') - \psi(T; x, u)\| &\leq T L_{\phi,u} e^{T L_{\phi,x}}\| u' - u \|. \\
    \end{aligned}
    \end{equation*}
    Finally, we pick $T_1 \vcentcolon= \min\{T', T''\}$ 
    and define $L_{\psi,x}\vcentcolon=e^{L_{\phi,x}T_1}L_{\phi,x}$ and 
    $L_{\psi_,u}\vcentcolon=e^{L_{\phi,x}T_1}L_{\phi,u}$, such that
    \begin{equation}
        \| \psi(T; x, u) - x\| \leq T \cdot (L_{\psi,x}\| x \| + L_{\psi,u} \| u \|)   
    \end{equation}
    and
    \begin{equation*}
    \begin{aligned}
        \|\psi(T; x, u') - \psi(T; x, u)\| \leq T  L_{\psi,u} \| u' - u\|.
    \end{aligned}
    \end{equation*}
    for any $x \in X_{\bar{V}}$, all $z,z' \in \mathcal{B}(\bar{z}(x), r_z)$ and any $T \leq T_1$.
    \hfill$\qed$

\section{Proof of Theorem 3.2.25}\label{app:ps_stab}
    It suffices to show that $V_l(w)$ is a Lyapunov function for \eqref{equ:ps} in $\mathbb{R}_{\geq0}^n$.
    Define $\hat{w}_{\min} \vcentcolon = \underset{i=1,\dots,n}{\min} \, [\,\hat{w}\,]_i$ and 
    $\hat{w}_{\max} \vcentcolon = \underset{i=1,\dots,n}{\max} \, [\,\hat{w}\,]_i$. The following inequalities hold:
    \begin{equation*}
        \hat{w}^{\top}w \geq \hat{w}_{\min} \cdot \mathbf{1}^{\top}w =
        \hat{w}_{\min} \cdot \| w \|_1 \geq \hat{w}_{\min} \cdot \| w \|_2
    \end{equation*}
    and 
    \begin{equation*}
        \begin{aligned}
            &&\hat{w}^{\top}w &\leq \hat{w}_{\max} \cdot \mathbf{1}^{\top}w \\
            && &\leq \hat{w}_{\max} \cdot \| \mathbf{1} \|_2  \| w \|_2 \\
            && &\leq \sqrt{n} \,\hat{w}_{\max}\cdot \| w \|_2,
    \end{aligned}
    \end{equation*}
    which show that there exist $\mathcal{K}_{\infty}$ functions 
    $\alpha_{l,1}(\| w \|) \vcentcolon= \hat{w}_{\min} \cdot \| w \|$ 
    and
    $\alpha_{l,2}(\| w \|) \vcentcolon= \sqrt{n}\, \hat{w}_{\max} \cdot \| w \|$ 
    such that
    \begin{equation}
        \alpha_{l,1}(\| w \|) \leq V_l(w) \leq \alpha_{l,2}(\|w\|).
    \end{equation}
    Moreover, for any $w>0$, we have that
    \begin{equation}
    \begin{aligned}
        &&V_l(w_+) - V_l(w) &= \hat{w}^{\top}Aw - \hat{w}^{\top}w \\
        &&                  &= \hat{w}^{\top}(A-\mathbb{I})\,w \\
        &&                  &= w^{\top}(A^{\top} - \mathbb{I})\,\hat{w} \\
        &&                  &\leq - \hat{d} \cdot \| w \|_1 \leq -\hat{d} \cdot \| w \|_2.
    \end{aligned}
    \end{equation}
    Hence, there exists a positive definite and continuous function 
    $\alpha_{l,3}(\| w \|)\vcentcolon=\hat{d} \cdot \| w \|$ such that, for any $w \geq 0$,
    the following holds
    \begin{equation}
        V_l(w_+) - V_l(w) \leq - \alpha_{l,3}(\| w \|)
    \end{equation}
    and $\alpha_{l,3}(0) = 0$, which concludes the proof. \hfill$\qed$
\begin{rem}
    Notice that the original Theorem in \cite{Kaczorek2008} requires the 
    existence of a strictly positive vector $\hat{w} > 0$ such that 
    \begin{equation}
        (A - \mathbb{I})\hat{w} < 0.
    \end{equation}
    Although the condition used in Theorem \ref{thm:ps_stab} is equivalent to one above, 
    the resulting $\hat{w}$ can only be used to define a Lyapunov function for the 
    dual system $w_+ = A^{\top}w$.
\end{rem}

\end{document}